\documentclass{amsart}
\usepackage{latexsym}
\newtheorem{theorem}{Theorem}[section]
\newtheorem{lemma}[theorem]{Lemma}
\newtheorem{proposition}[theorem]{Proposition}
\newtheorem{corollary}[theorem]{Corollary}
\theoremstyle{definition}                               
\newtheorem{definition}[theorem]{Definition}
\newtheorem*{example*}{Example}
\newtheorem{example}[theorem]{Example}

\theoremstyle{remark} 
\newtheorem{remark}[theorem]{Remark}
\newtheorem*{remark*}{Remark}
\numberwithin{equation}{section}
\DeclareMathOperator{\Grad}{\mbox{\sf Grad}}
\DeclareMathOperator{\Mad}{\mbox{\sf Mad}}
\DeclareMathOperator{\Fad}{\mbox{\sf Fad}}

\DeclareMathOperator{\Pic}{\sf Pic}

\DeclareMathOperator{\End}{\sf End}
\DeclareMathOperator{\Aut}{\sf Aut}
\DeclareMathOperator{\Spec}{\sf Spec}
\DeclareMathOperator{\ad}{ad}

\DeclareMathOperator{\Frac}{\sf Frac}
\DeclareMathOperator{\gr}{gr}
\DeclareMathOperator{\Ker}{Ker}
\def\Gr{\mbox{\rm{Gr}}^{\mbox{\scriptsize{\rm{ad}}}}}
\def\Grat{\mbox{\rm{Gr}}^{\mbox{\scriptsize{\rm{rat}}}}}
\def\O{\mathcal{O}}
\def\K{\mathbb{K}}

\def\C{\mathcal{C}}
\def\P{\mathbb{P}}
\def\A{\mathbb{A}}
\def\B{\mathcal{B}}
\def\L{\mathcal{L}}
\def\D{\mathcal{D}}
\def\c{\mathbb{C}}
\def\d{\partial}
\begin{document}
\title[Mad subalgebras]
{Mad subalgebras of rings of differential operators on curves}
%
%
	\author{Yuri Berest}
%
%
	\address{Department of Mathematics, Cornell University,
         Ithaca, NY 14853, USA}
	\email{berest@math.cornell.edu}
%
     \author{George Wilson}
     \address{Mathematical Institute, 24--29 St Giles, Oxford OX1 3LB, UK}
      \email{wilsong@maths.ox.ac.uk}
%
%
\begin{abstract}
We study the maximal abelian ad-nilpotent (mad) subalgebras of the 
domains $\, \D $ Morita equivalent to the first 
Weyl algebra. We give a complete 
description both of the individual mad subalgebras and of the space of 
all such. A surprising consequence is that this last space is 
independent of $\, \D \,$. Our results generalize some classic theorems of 
Dixmier about the Weyl algebra.
\end{abstract}
\maketitle
\section{Introduction and statement of results}

We begin by recalling some results of Dixmier (see \cite{D}) about the (first) 
Weyl algebra $\, A \,$.  We shall think of $\, A \,$ as the algebra 
$\, \D(\A^1) \,$ of differential operators on the (complex) affine line, 
that is, 
as the algebra $\, \c[z, \d_z] \,$ of polynomial differential operators in 
one variable $\, z \,$. We call an element $\, b \in A \,$ {\it ad-nilpotent} 
(``strictly nilpotent'' in \cite{D}) if for each $\, a \in A \,$ we have 
$\, (\ad b)^k(a) = 0 \text{\, for some \,} k \,$.  We call 
a maximal abelian subalgebra $\, B \,$ of  $\, A \,$ a {\it mad subalgebra} if 
every element of $\, B \,$ is ad-nilpotent.  For example, $\, \c[z] \,$ is 
clearly a mad subalgebra of $\, A \,$, and $\, \c[\d_z] \,$ is another. 
One of Dixmier's main aims in \cite{D} was to 
obtain information about the group $\, \Aut A \,$ of 
$\, \c$-automorphisms of $\, A \,$; to this end he studied the action of 
$\, \Aut A \,$ on the set of mad subalgebras of $\, A \,$. One of his key  
results was that this action is {\it transitive}.  Clearly, that implies
\begin{theorem}[\cite{D}]
\label{D1}
Every mad subalgebra $\, B \subset A \,$ has the form 
$\, B = \c[x] \,$ for some $\, x \in A \,$.
\end{theorem}

If $\, B \,$ is a mad subalgebra, we shall call a choice of generator 
$\, x \,$ for $\, B \,$ a {\it framing} of $\, B \,$, and the pair 
$\, (B,x) \,$ a {\it framed mad subalgebra} of $\, A \,$.  Dixmier 
showed in fact (see \cite{D}, Lemme 8.9) 
that $\, \Aut A \,$ acts transitively on the set 
$\, \Mad A \,$ of all framed mad subalgebras of $\, A \,$. Let 
$\, \Gamma \,$ be the subgroup of $\, \Aut A \,$ consisting of all 
automorphisms $\, \gamma_p \,$ of the form
\begin{equation}
\label{gamma}
\gamma_p(z) = z \,, \quad \gamma_p(\d_z) = \d_z - p'(z)
\end{equation} 
where $\, p \in \c[z] \,$ (we may think of $\, \gamma_p \,$  as 
conjugation by $\, e^{p(z)} \,$).  It is easy to check that $\, \Gamma \,$ 
is exactly the isotropy group of $\, z \in A \,$, or, equivalently, 
of the natural base-point $\, (\c[z], z) \in \Mad A \,$.  Dixmier's result 
can therefore be formulated as follows. 
\begin{theorem}[\cite{D}]
\label{D2}
There is a natural bijection 
$$
\Aut(A)/ \Gamma \to \Mad A \,.
$$
\end{theorem}

We now wish to generalize Theorems~\ref{D1} and \ref{D2} to the case where 
the Weyl algebra is replaced by the ring $\, \D(X) \,$ of differential 
operators on any affine curve.  Clearly, the ring 
$\, \O(X) \,$ of regular functions on $\, X \,$ is a mad subalgebra 
of $\, \D(X) \,$.  Theorems of Makar-Limanov and Perkins (see \cite{M}, 
\cite{P}) show that $\, \O(X) \,$ is the {\it only} mad subalgebra 
except in the case when $\, X \,$ is a {\it framed curve}, by which we mean 
that there is a regular bijective map $\, \pi : \A^1 \to X \,$  
(thus topologically a framed curve is simply the affine line,
but it may have an arbitrary finite number of cusps).  From now on 
we suppose that  $\, X \,$ is a framed curve, since 
Theorems~\ref{D1} and \ref{D2} can have interesting generalizations 
only in that case.  Our first result is as follows.
\begin{theorem}
\label{T1}
Let $\, B \,$ be any mad subalgebra of $\, \D(X) \,$, where $\, X \,$ 
is a framed curve. Then $\, \Spec B \,$ is a framed curve.
\end{theorem}

This Theorem is a sharper version of a result of \cite{LM}, where it 
is shown that the normalization of $\, \Spec B \,$ is always isomorphic 
to $\, \A^1 $; that means that  $\, \Spec B \,$ is obtained from 
$\, \A^1 $ by introducing cusps, but also, perhaps, identifying certain 
points of $\, \A^1 $  to form double points, 
or even higher order multiple points.  The new part of Theorem~\ref{T1}
is thus the assertion that multiple points do not occur.
The question of whether $\, \Spec B \,$ is necessarily 
a framed curve (that is, free of multiple points) 
was raised by P.~Perkins (see \cite{P}) 
in a special case where the mad subalgebra $\, B \,$ is dual 
(in the sense of Section~\ref{dual} below) to $\, \O(X) \,$. 
He raised also a more subtle question: setting $\, Y := \Spec B \,$, is it 
true that $\, \D(X) \,$ is isomorphic to $\, \D(Y) \,$?  
Our proof of Theorem~\ref{T1} yields also the answer to this question, 
namely ``not always''; more precisely:
\begin{theorem}
\label{T2}
Let $\, B \,$ be any mad subalgebra of $\, \D(X) \,$, and let 
$\, Y := \Spec B \,$.  Then there is a rank {\rm 1} torsion-free 
coherent sheaf $\, \L \,$ over $\, Y $ and an isomorphism 
$\, \varphi : \D_{\L}(Y) \to \D(X) \,$ such that 
$\, \varphi(\O(Y)) = B \,$.  
\end{theorem}

Here $\, \D_{\L}(Y) \,$ denotes the ring of differential operators 
on global sections of $\, \L \,$. If $\, \L \,$ is not locally free, 
then $\, \D_{\L}(Y) \,$ is 
not necessarily isomorphic to $\, \D(Y) \,$ (see \cite{BW3}, Example 8.4).

Theorems \ref{T1} and \ref{T2} give a satisfactory description of 
the individual mad subalgebras $\, B \subset \D(X) \,$; we now describe 
the ``space'' of all such $\, B \,$, in the spirit of 
Theorem~\ref{D2}. As we saw above, if $\, B \,$ is any mad subalgebra 
of $\, \D(X) \,$, then its integral closure $\, \overline{B} \,$ 
is isomorphic to a 
polynomial algebra $\, \c[x] \,$; as before, we call a choice of 
generator for $\, \overline{B} \,$ a {\it framing} of $\, B \,$, 
and we denote the set of all framed mad subalgebras of $\, \D(X) \,$ 
by $\, \Mad \D(X) \,$. Generalizing Theorem~\ref{D2}, we shall prove
\begin{theorem}
\label{T3}
For any framed curve $\, X \,$, there is a natural bijection 
$$
\Aut(A) / \Gamma \to \Mad \D(X) \,.
$$
\end{theorem}

We found this result surprising: it implies that the space of mad subalgebras 
of $\, \D(X) \,$ is independent of $\, X \,$. The algebras 
$\, \D(X) \,$ are (up to isomorphism) 
exactly the domains Morita equivalent to the 
Weyl algebra; thus $\, \Mad \D \,$ is a Morita invariant 
for this special class of algebras.  
It would be interesting to understand whether this is an instance of some 
more general principle.

The last theorem that we want to formulate in this Introduction 
describes the quotient space of 
$\, \Mad \D(X) \,$ by the natural action of the automorphism group 
of $\, \D(X) \,$.  Recall (see \cite{W2}) that for each $\, n \geq 0 \,$ the 
{\it Calogero-Moser space} $\, \C_n \,$ is the space of isomorphism 
classes of triples $\, (\mathcal{V}; \mathbb{X},\mathbb{Y}) \,$ 
where $\, \mathcal{V} \,$ is an 
$n$-dimensional complex vector space, 
and $\, \mathbb{X} \,$ and $\, \mathbb{Y} \,$ 
are endomorphisms of $\, \mathcal{V} \,$ 
such that $\, [\mathbb{X},\mathbb{Y}] + I \,$ has rank 1\,. 
We make the group $\, \Gamma \,$ act on $\, \C_n \,$ by the formula
\begin{equation}
\label{Caction}
\gamma_p(\mathbb{X},\mathbb{Y}) = (\mathbb{X} + p'(\mathbb{Y}), \mathbb{Y}) \,.
\end{equation}
Recall further (see \cite{K1}, \cite{K2}, \cite{BW1}, \cite{BW3}) 
that the algebras $\, \D(X) \,$ 
are classified up to isomorphism by a non-negative integer $\, n \,$ which 
we call the {\it differential genus} of $\, X \,$: it can be 
thought of as the number of cusps of $\, X \,$, but 
counted with appropriate multiplicities (see formula (8.3) in \cite{BW3}).

\begin{theorem}
\label{T4}
Let $\, X \,$ be a framed curve, and let $\, n \,$ be its differential 
genus.  Then there is a natural bijection
$$
\C_n / \Gamma \to {\Mad \D(X)} / {\Aut \D(X)} \,.
$$
\end{theorem}

This theorem was announced (without proof) in \cite{BW3}.  
The space $\, \C_0 \,$ is a point, so in the case where $\, X = \A^1 \,$ 
Theorem~\ref{T4} reduces to Dixmier's result that $\, \Aut A \,$ acts 
transitively on $\, \Mad A \,$.  In general, $\, \C_n \,$ is a smooth 
affine variety of dimension $\, 2n \,$, and generic orbits of $\, \Gamma \,$ 
are $n$-dimensional, so the theorem suggests that 
$\, {\Mad \D(X)} / {\Aut \D(X)} \,$ 
is an $n$-dimensional space.  Unfortunately, 
we do not know any intrinsic 
way of assigning a dimension to this space; and even  
$\, \C_n / \Gamma \,$ is not a good quotient in the sense of algebraic 
geometry (for example, because $\, \Gamma \,$ has some orbits of dimension 
less than $\, n \,$, at least for $\, n>2 \,$).

Despite its modest appearance, Theorem~\ref{T1} is the key result of this 
paper, the others being comparatively formal consequences. 
Its proof involves a curious mixture of 
familiar algebraic arguments and others connected with the theory of
integrable systems; in particular, the Burchnall-Chaundy theory of 
commuting ordinary differential operators plays a crucial role.  
We offer two versions of the proof, in one of which (playing devil's 
advocate) we have sought to reduce the role of the Burchnall-Chaundy theory 
to a minimum. We do not know how to eliminate it entirely: 
we leave the possibility of that as a worry for the reader. A detailed 
overview of the contents of the paper can be found in the introductory 
remarks to the the individual sections that follow.  Here we just mention 
that Sections~\ref{diffop} and \ref{several} make no claim to 
originality, but are an exposition of some of the results of \cite{LM}. 
We give the exposition in some detail, because we rely heavily on these results, 
and the account of them given in \cite{LM} is not quite  
satisfying (specifically, part of Section 3 of that paper is missing, and 
the reference to \cite{P} in the proof of its Corollary 4.6 seems 
difficult to justify directly). The version presented here is based on 
notes graciously placed at our disposition by G.\ Letzter.
\\*[.5ex]
\noindent
\scriptsize 
\textbf{Acknowledgments}. We thank G.\ Letzter for a very helpful 
correspondence concerning the material in Section~\ref{several} below,  
and also for pointing out some errors in the  first version of this paper. 
We are indebted to the referee for his careful reading of the paper, 
and for his thoughtful suggestions for improving the exposition.    
The authors were partially supported by 
NSF grant DMS 04-07502 and an A.\ P.\ 
Sloan Research Fellowship; the second author is grateful to the 
Mathematics Department of Cornell University for its hospitality 
during the preparation of this article.   
\normalsize

\section{Mad subalgebras}
\label{defmad}

In this section we give some definitions, including those of 
mad subalgebras and filtrations of an algebra $\, A \,$: 
these abstract the basic properties 
of the rings $\, \D(X) \,$ of differential operators on algebraic 
varieties (see Example~\ref{dop} below).  We also note some special 
features of the ``$1$-dimensional'' case where $\, A \,$ satisfies 
the condition \eqref{cc} below.

Let $\, A \,$ be a noncommutative algebra over $\, \c \,$.  As usual, 
for each $\, b \in A \,$ we write $\, \ad b \,$ for the 
inner derivation of $\, A \,$ defined by $\, (\ad b)(a) := [b,a] \,$; 
we set
$$
N_k(b) := \Ker(\ad b)^{k+1}\ ,\ \ N(b) := \bigcup_{k \geq 0} N_k(b) \,.
$$ 
It is easy to check that $\, N(b) \,$ is a filtered subalgebra of 
$\, A \,$.  If $\, N(b) = A \,$, we say that $\, b \,$ is a 
{\it (locally) ad-nilpotent} element of $\, A \,$, and we call the 
above filtration on $\, A = N(b) \,$ the {\it filtration induced 
on $\, A \,$ by $\, b \,$}. In later sections we shall sometimes write 
$\, N_A(b) \,$ instead of $\, N(b) \,$ (if the algebra $\, A \,$ is not 
clear from the context).
For $\, b \in A \,$, we denote the centralizer of $\, b \,$ by $\, C(b) \,$ 
(or, if necessary, by $\, C_A(b) \,$); 
thus $\, C(b) \equiv N_0(b) \,$ as defined above.
\begin{proposition}
\label{CN}
Suppose $\, C(b_1) = C(b_2) \,$.  Then $\, N(b_1) = N(b_2) \,$ 
as filtered algebras.
\end{proposition}

The proof depends on the following (purely set-theoretical) lemma.
\begin{lemma}
\label{commder}
Let $\, f_1 , f_2 : A \to A\,$ be two maps 
such that (i) $\, f_1 \,$ and $\, f_2 \,$ commute; 
(ii) $\, \Ker f_1 = \Ker f_2 \,$. 
Then $\, \Ker f_1^n = \Ker f_2^n \,$ for all 
$\, n \geq 1 \,$.
\end{lemma}
\begin{proof}
An easy induction on $\, n \,$ ($ f_1 \,$ and $\, f_2 \,$ do not 
even have to be linear). 
\end{proof}
\begin{proof}[Proof of Proposition~\ref{CN}]
The elements $\, b_1 \,$ and $\, b_2 \,$ commute, hence the derivations 
$\, \ad b_1 \,$ and $\, \ad b_2 \,$ commute.  So Lemma~\ref{commder} 
applies to give  $\, N_k(b_1) = N_k(b_2)$ for all $\, k \geq 0 \,$, which is 
what the Proposition asserts.
\end{proof}

More generally, if $\, B \,$ is any subset of $\, A \,$, we can define 
the filtered subalgebra 
$\, N(B) = \bigcup_{k \geq 0} N_k(B) \,$, where
$$
N_k(B) := \{ a \in A : (\ad b_0) (\ad b_1) \ldots (\ad b_k)(a) = 0 \ \ 
\text{for all} \ 
\, b_0, b_1 \ldots b_k \in B \} \,.
$$
We are interested in the case when $\, B \,$ is an abelian subalgebra 
of $\, A \,$: we say $\, B \,$ is {\it ad-nilpotent} if 
$\, N(B) = A \,$.   Choosing $\, b_0 = b_1 = \ldots = b_k \,$ in 
the definition of $\, N_k(B) \,$, we see that 
if $\, B \,$ is ad-nilpotent, then every element of $\, B \,$ is ad-nilpotent. 
If $\, B \,$ is ad-nilpotent, we call the natural filtration on 
$\, N(B) \equiv A \,$ the {\it filtration induced by $\, B \,$}.
Clearly, in this filtration $\, \{ A_k \} \,$ the ring $\, A_0 \,$ 
is the {\it commutant} $\, C(B) \,$ of $\, B \,$.
\begin{definition}
\label{mad}
We say that $\, B \,$ is a {\it mad} subalgebra of $\, A \,$ if \\
(i)  $\, B \,$ is ad-nilpotent; \\
(ii)  $ C(B) = B \,$ .
\end{definition} 

If  $\, B \,$ is a mad subalgebra of $\, A \,$, then the filtration 
$\, \{ A_k \} \,$ induced by $\, B \,$ has the properties
\begin{enumerate}
\item{$ \, A_{-1} = 0 \,$ (that is, the filtration is {\it positive});}
\item{if $\, b \in A_0, \, a \in A_k \,$, then $\, [b,a] \in A_{k-1} \,$;}
\item{if $\, a \,$ has filtration degree $\, k \,$, then there is a 
$\, b \in A_0 \,$ such that $\, [b,a] \,$ }has filtration degree $\, k-1 \,$.
\end{enumerate}
We call a filtration of $\, A \,$ with these properties a 
{\it mad filtration}.  It is easy to check that if $\, \{ A_k \} \,$ 
is a mad filtration and we define $\, B := A_0 \,$, then 
$\, B \,$ is a mad subalgebra of $\, A \,$, and the given filtration 
coincides with the one induced by $\, B \,$.  In this way the mad 
subalgebras of $\, A \,$ correspond 1--1 to the mad filtrations.
\begin{example}
\label{dop}
Our definitions of mad subalgebras and filtrations are modelled on 
the following situation.  Let $\, X \,$ be an irreducible complex 
affine variety, $\, \O \,$ the ring of regular functions on  $\, X \,$, 
and let $\, E := \End_{\c} \O \,$.  The filtered algebra 
$\, N_E(\O) \,$ is (by definition) the ring $\, \D(X) \,$ of 
{\it differential operators} on $\, X \,$, and $\, \O \,$ is a  
mad subalgebra of $\, \D(X) \,$.  More generally, let $\, \L \,$ 
be a rank 1 torsion-free coherent sheaf over $\, X \,$, and 
$\, M \,$ the corresponding $\, \O $-module (of global sections); 
if $\, E := \End_{\c} M \,$, then the filtered algebra 
$\, N_E(\O) \,$ is (by definition) the ring $\, \D_{\L}(X) \,$ of 
differential operators on $\, \L \,$.  Here it may happen that 
the centralizer $\, C_E(\O) \,$ is slightly larger than $\, \O \,$: 
we call $\, \L \,$ {\it maximal} if $\, C_E(\O) = \O \,$.  In that 
case  $\, \O \,$ is a mad subalgebra of $\, \D_{\L}(X) \,$.  
\end{example}
\begin{remark*}
Every line bundle (locally free rank 1 coherent sheaf) over $\, X \,$ 
is maximal, but the converse is not true (see, for example \cite{SW}, 
p.\ 46).  For this reason our notion of a ``ring with mad filtration'' is 
slightly more general than the ``algebras of twisted differential 
operators'' introduced in \cite{BB} (which model the case where 
$\, \L \,$ is a line bundle).
\end{remark*}
We suppose from now on that our algebra $\, A \,$ satisfies the condition
\begin{equation}
\label{cc}
C(a) \text{ is commutative for each } \, a \in A \setminus \c \,.
\end{equation}
This condition is very restrictive; for example, if $\, A \,$ is the 
ring of differential operators on an affine variety $\, X \,$, then 
\eqref{cc} is satisfied only if $\, X \,$ is 1-dimensional. However, that is 
exactly the situation that concerns us in this paper. 
Many things become simpler if \eqref{cc} holds: for example, 
we have $\, B = C(B) \,$ if (and only if) 
$\, B \,$ is a maximal abelian subalgebra of $\, A \,$. Further, 
the maximal abelian subalgebras of $\, A \,$ are 
exactly the centralizers of the elements $\, b \in A \setminus \c \,$, and 
$\, C(b) \,$ is the unique maximal abelian subalgebra containing $\, b \,$. 
It follows that the intersection of any two distinct maximal abelian 
subalgebras is $\, \c \,$.  The following facts about {\it mad} 
subalgebras are all easy consequences of 
\eqref{cc} and Proposition~\ref{CN}. 
\begin{proposition}
Suppose that $\, A \,$ satisfies \eqref{cc}, and let 
$\, b \in A \setminus \c \,$ be ad-nilpotent. Then $\, C(b) \,$ is 
a mad subalgebra of $\, A \,$.
\end{proposition}
\begin{proposition}
Suppose that $\, A \,$ satisfies \eqref{cc}, and let $\, B \,$ be a 
mad subalgebra of $\, A  \,$.  Then  the filtration induced on 
$\, A \,$ by $\, B \,$ coincides with the filtration induced by 
any element of  $\, B \setminus \c \,$. In particular, if $\, a \,$ has 
degree $\, k > 0 \,$ in this filtration, then $\, [b,a] \,$ has 
degree (exactly) $\, k-1 \,$ for \emph{every} 
$\, b \in B \setminus \c \,$.
\end{proposition}
\begin{proposition}
Suppose that $\, A \,$ satisfies \eqref{cc}, and let $\, B \,$ be maximal 
among the abelian subalgebras of $\, A \,$ all of whose elements are 
ad-nilpotent.  Then $\, B \,$ is a mad subalgebra of $\, A \,$.
\end{proposition}

These propositions indicate that if \eqref{cc} is satisfied, then 
various possible definitions of ``mad'' all coincide; in particular, 
the definition given in the present section agrees with the one we 
used in the Introduction.
\section{Rings of Differential Operators}
\label{diffop}

The next two sections provide a self-contained exposition of some of the 
results of \cite{LM}.  The present section gathers together some 
preliminary facts, culled from \cite{S}, \cite{M}, \cite{KM} and 
\cite{LM}. The main points to note are Proposition~\ref{QC}, which 
ensures that the algebras studied later on all satisfy the 
condition \eqref{cc}; and the more technical-looking 
Proposition~\ref{adbar}, which lies at the heart of the proofs in 
the following Section~\ref{several}. 

Let $\, \K \,$ be a commutative field 
containing\footnote{For the purposes of this section we could consider 
that $\c$ denotes an arbitrary field of characteristic zero.} 
$\, \c \,$, and let $\, \d \,$ be a derivation of $\, \K \,$ 
with kernel $\, \c \,$.  Then we can form the ring 
$\, \K[\d] \,$, consisting of expressions of the form 
$$
D = \sum_0^n f_i \d^i \ , \ \ f_i \in \K \,,
$$
with multiplication defined by the commutation relation
$$
[\d,f] =\d(f) \text{ \ for all \,} f \in \K \,.
$$
Clearly, the ring $\, \K[\d] \,$ does not change if we replace 
$\, \d \,$ by $\, f\d \,$ for some non-zero $\, f \in \K \,$.
We have in mind principally the case when $\, \K \,$ is the 
function field of a curve, so that $\, \K \,$ is an extension of 
$\, \c \,$ of transcendence degree $\, 1 \,$.  In that case 
the $ \c$-derivations of $\, \K \,$ form a $1$-dimensional 
$\, \K $-vector space, so the algebra $\, \K[\d] \,$ has an 
intrinsic interpretation (independent of the choice of $\, \d \,$) 
as the ring $\, \D(\K) \,$ of differential operators on $\, \K \,$. 
 
It is easy to show that $\, \K[\d] \,$ is a Noetherian domain, 
hence it has a 
quotient field $\, Q \,$.  
It is sometimes helpful to think of $\, Q \,$ 
as sitting inside the still larger field 
$\, \overline{Q} = \K((\d^{-1})) \,$ of formal 
Laurent series
\begin{equation}
\label{laurent}
D = \sum_{- \infty}^n f_i \d^i \ , \ \ f_i \in \K \,.
\end{equation}
If $\, D \,$ has the form \eqref{laurent} with $\, f_n \not= 0 \,$, we call 
$\, f_n \d^n \,$ the {\it leading term} of $\, D \,$ and $\, f_n \,$ its 
{\it leading coefficient}.  The following fact goes back to Schur 
(see \cite{S}).
\begin{proposition}
\label{QC}
$\, \overline{Q} \,$ satisfies the condition \eqref{cc}.
\end{proposition}
\begin{proof}[Sketch of proof]
There are three cases.\\
(i) Suppose that $\, L \in \overline{Q} \,$ has leading term $\, a \d^n \,$, 
where $\, n \not= 0 \,$.  If necessary we adjoin to $\, \K \,$ an $n$th 
root $\, \alpha \,$ of $\, a \,$.  Then one shows that $\, L \,$ has an 
$n$th root $\, L^{1/n} = \alpha \d + \ldots \,$ in 
$\, \K(\alpha)((\d^{-1})) \,$, and that the centralizer of $\, L \,$ 
in this field
consists of the Laurent series in $\, L^{1/n} \,$.  Clearly, this 
is commutative. For more details, see \cite{S}. \\
(ii) Suppose $\, L = a + a_1 \d^{-1} + \ldots \,$ has order $\, 0 \,$ with 
$\, a \in \c \,$ (and $\, L \not= a \,$).  Then 
$\, C(L) = C(L-a) \,$, which is commutative by case (i).\\
(iii)  Suppose $\, L = a + a_1 \d^{-1} + \ldots \,$ has order $0$ with 
$\, a \notin \c \,$.  If $\, P \,$ has leading term $\, p \d^m $ 
with $\, m \not= 0 \,$, 
then $\, [P,L] \,$ has leading coefficient 
$\, mp \d(a) \not= 0 \,$; hence $\, C(L) \,$ consists of operators 
of order zero. If now 
$\, P_1, P_2 \in C(L) \,$, then $\, [P_1, P_2] \in C(L) \,$ is either 
zero or an operator of order $\, < 0 \,$. We just saw that the latter is 
impossible, hence $\, C(L) \,$ is commutative.  Alternatively, to 
get a more precise result, we can argue as follows: equating 
coefficients of powers of $\, \d \,$ in the expansion of $\, PL = LP \,$ 
shows that for each $\, p \in \K \,$ there is a unique 
operator of the form $\, P = p + p_1 \d^{-1} + \ldots \,$ that commutes 
with $\, L \,$. It follows that  $\, C(L) \,$ is isomorphic 
to $\, \K \,$.
\end{proof}

Of course, it follows from Proposition~\ref{QC} 
that any subalgebra of $\, \overline{Q} \,$, in particular $\, Q \,$, 
satisfies \eqref{cc}.

Our aim in the rest of this section is to prove the basic 
Proposition~\ref{adbar} below.  We start with 
the following very simple lemma.
\begin{lemma}
\label{q}
Let $\, F \,$ be a field of characteristic $\, 0 \,$ (not necessarily 
commutative), and let $\, \d \,$ be a derivation of $\, F \,$.  Suppose that 
$\, \Ker \, \d^2 \not= \Ker \, \d \,$. Then there is a $\, q \in F \,$ such 
that $\, \d(q) = 1 \,$; and for each $\, n \geq 1 \,$, $\, \Ker \, \d^n \,$ 
is an $n$-dimensional (left or right) vector space over  $\, \Ker \, \d \,$ 
with basis $\, \{ 1, q, \ldots, q^{n-1} \} \,$. 
\end{lemma}
\begin{proof}
Choose $\, r \in \Ker \, \d^2 \setminus \Ker \, \d \,$, and let 
$\, s = \d(r) \,$, so that $\, s \not= 0 \,$ but $\, \d(s) = 0 \,$.  
Set $\, q = s^{-1}r \,$; then  $\, \d(q) = 1 \,$.  The rest is an 
easy induction on $\, n \,$ 
(using the fact that $\, \d(q^n) = nq^{n-1} \,$).
\end{proof}

In particular, we can apply Lemma~\ref{q} in the case where 
$\, \d = \ad u \,$ for some $\, u \in F \,$; in this case it is 
tempting to denote the element $\, q \,$ in the Lemma by $\, -\d_u \,$.  
We record the result for future reference.
\begin{corollary}
\label{d}
Let $\, u \in F \,$ (where $\, F \,$ is a noncommutative field of 
characteristic $\, 0 $), and suppose that $\, N_F(u) \neq C_F(u) \,$.  
Then there is an element $\, \d_u \in F \,$ such that 
$\, [\d_u,u] = 1 \,$ and $\, N_F(u) = C_F(u)[\d_u] \,$.
\end{corollary}

We return now to our ring $\, \K[\d] \,$.  The next lemma is a special 
case of Proposition~\ref{adbar}.

\begin{lemma}
\label{ad}
Suppose that $\, D \in \overline{Q }\,$ has leading term 
$\, \d^n $, where $\, n \not= 0 \,$, 
and suppose $\, D \,$ acts ad-nilpotently on some operator $\, \Theta \,$ 
with leading term $\, f \d^m $, where $\, f \notin \c \,$.  Then the 
equation  $\, \d(q) = 1 \,$ has a solution in $\, \K \,$, and if $\, L \,$ 
is any operator on which $\, D \,$ acts ad-nilpotently, then the leading 
coefficient of $\, L \,$ belongs to $\, \c[q] \,$. 
\end{lemma}
\begin{proof}
We have $\, [\d^n, f \d^m] = n \d(f) \d^{n+m-1} + $ (lower order terms), 
hence for any $\, i \geq 1 \, $ the coefficient of $\, \d^{m+ i(n-1)} \,$ 
in $\, (\ad D )^i(\Theta) \,$ 
is $\, n^i \d^i(f) \,$.  So if  $\, D \,$ acts 
ad-nilpotently on $\, \Theta \,$, 
we have $\, \d^i(f) = 0 \,$ for some $\, i \geq 1 \,$, so that 
$\, \Ker \, \d^i \not= \c \,$.  Hence 
$\, \Ker \, \d^2 \not= \Ker \, \d \ (= \c) \,$, so  
Lemma~\ref{q} tells us that $\, q \,$ exists as stated, and that 
$\, f \in \c[q] \,$.  The last assertion in the Lemma is trivial 
if the leading coefficient of $\, L \,$ is a scalar, and otherwise 
follows by the argument above (applied to  $\, L \,$ instead 
of $\, \Theta \,$).
\end{proof}

Finally, we want to remove the hypothesis in Lemma~\ref{ad} 
that $\, D \,$ has scalar leading 
coefficient. This assumption is not essential, because 
we can always reduce to that case by a 
``change of variable''.  Recall that if $\, \widehat \K \,$ is 
a finite extension field of $\, \K \,$, then $\, \d \,$ extends 
uniquely to a derivation of $\, \widehat \K \,$, still with kernel 
$\, \c \,$: we denote this extension by the same symbol $\, \d \,$.
If $\, D \in \overline{Q}(\K) \,$ has 
leading term $\, a \d^n \,$, where $\, n \not= 0 \,$,
we can form the extension field $\, \widehat \K  = \K(\alpha) \,$, 
where $\, \alpha^n = a \,$.  Then $\, d := \alpha \d \,$ is a derivation of 
$\, \widehat \K \,$, and we may write the elements of 
$\, \overline{Q}(\widehat{\K}) \,$ as Laurent series in $\, d \,$ 
(rather than $\, \d \,$). The operator $\, D \,$ then has leading term 
$\, d^n $, so we may apply Lemma~\ref{ad} to $\, (\widehat \K, d) \,$ 
to get the following.
\begin{proposition}
\label{adbar}
Suppose that $\, D \in \overline{Q} \,$ 
has leading term $\, a \d^n \,$, where $\, n \not= 0 \,$, 
and suppose $\, D \,$ acts ad-nilpotently on some operator 
$\, \Theta \in \overline{Q} \,$ 
with leading term $\, f \d^m $, where $\, f^n/a^m \notin \c \,$.  
Let $\, \widehat \K = \K(\alpha) \,$, where $\, \alpha^n = a \,$. 
Then the equation  $\, \alpha \d(q) = 1 \,$ has a solution  
$\, q \in \widehat{\K} \,$, and if $\, L \,$ 
is an operator with leading term $\, \beta \d^r $
on which $\, D \,$ acts ad-nilpotently, 
then $\, \beta \in \alpha^r \c[q] \,$. 
\end{proposition}

\section{Rings with several mad subalgebras}
\label{several}

In this section we conclude our reworking of some parts of 
\cite{LM} which were not treated convincingly in that paper. 
The main results are Theorems~\ref{x} and \ref{y}.  

For the rest of the paper $\, Q \,$ will denote 
the quotient field of the Weyl 
algebra, and $\, \D \,$ will be a subalgebra of $\, Q \,$ with the 
properties 
\begin{equation}
\label{ass1}
\text{the quotient field of} \ \, \D \, \ \text{is} \  \, Q \,;
\end{equation}
\begin{equation}
\label{ass2}
\, \D \,\  \text{contains more than one mad subalgebra}.
\end{equation}
We fix a mad subalgebra $\, B \subset \D \,$. 
We may regard its field of fractions 
$\, \Frac B \,$ as a subfield of $\, Q \,$; in particular, the 
integral closure $\, \overline{B} \,$ of $\, B \,$ (in $\, \Frac B \,$) 
is a subalgebra of $\, Q \,$.
\begin{theorem}
\label{x}
There is an $\, x \in Q \,$ such that $\, \overline{B} = \c[x] \,$.
\end{theorem}

The proof uses the following lemma, which is well known (see, 
for example \cite{E2}, \cite{J} p.\ 256, \cite{P} p.\ 281). 
However, we shall give a self-contained proof.
\begin{lemma}
\label{pol}
Let $\, B \not= \c \,$ be a subalgebra of a polynomial algebra 
$\, \c[q] \,$. Then 
(i)  $\, B \,$ is finitely generated; (ii) the integral closure 
$\, \overline{B} \,$ of 
$\, B \,$ has the form $\, \c[x] \,$ for some $\, x \in \c[q] \,$.
In other words, $\, B \,$ is the coordinate ring of a curve with 
normalization isomorphic to $\, \A^1 \,$.
\end{lemma}
\begin{proof}
(i) follows from the fact that 
every sub-semigroup of $\, \mathbb{N} \,$ is finitely generated (the 
degrees of the polynomials in $\, B \,$ form such a semigroup), while 
we can see (ii) from general principles as follows.  By L\"uroth's 
Theorem, $\, \Spec B \,$ is a rational curve, hence 
$\, \Spec \overline{B} \,$ is isomorphic to  $\, \A^1 \,$ with 
(perhaps) a finite number of points removed.  Because $\, \c[q] \,$ is 
integrally closed, we have 
$\, \overline{B} \subseteq \c[q] \,$: this inclusion corresponds to a map 
$\, f: \A^1 \to \Spec \overline{B} \,$ with dense image 
$\, \A^1 \setminus S \,$ for some finite set $\, S \,$. 
To see that  $\, S \,$ is empty, we regard $\, f \,$ as a 
rational map $\, \P^1 \to \P^1 \,$. Because  $\, \P^1 \,$ is 
a smooth curve, this map is  
regular everywhere; so 
its image is closed, hence equal to $\, \P^1 \,$; and the point 
at infinity maps to only one point.  Thus $\, f \,$ maps $\, \A^1 \,$ 
onto $\, \A^1 \,$; that is, $\, S \,$ is empty and 
$\, f(\A^1) = \Spec \overline{B}  \simeq \A^1 \,$.
\end{proof}
\begin{proof}[Proof of Theorem~\ref{x}]
Let $\, \K \,$  be the centralizer of $\, B \,$ in $\, Q \,$; by 
Proposition~\ref{QC}, $\, \K \,$ is a commutative field.  Choose 
any $\, u \in \K \setminus \c \,$; then $\, C_Q(u) = \K \,$ is 
commutative and $\, N_Q(u) = N_Q(B) \,$ is not 
(because it contains $\, \D \,$, 
which is not commutative by \eqref{ass1}).  So by Corollary~\ref{d} 
we can choose $\, \d_u \in Q \,$ such that $\, [\d_u,u] = 1 \,$ and 
$\, N_Q(B) = \K[\d_u] \,$.  The derivation 
$\, \ad \d_u \,$ of $\, \K \,$ has kernel $\, \c \,$, for this kernel 
is the intersection of $\, C_Q(u) \,$ and $\, C_Q(\d_u) \,$; 
since $\, u \,$ and 
$\, \d_u \,$ do not commute, their centralizers are distinct, and hence 
have intersection $\, \c \,$.  We may therefore 
think of $\, Q \,$ as embedded in the field 
$\, \K((\d_u^{-1})) \,$ and apply the results of Section~\ref{diffop}.

By assumption \eqref{ass2}, we may choose an ad-nilpotent element 
$\, D \in \D \setminus B \,$. Using the inclusion 
$\, \D \subseteq \K[\d_u] \,$, we think of $\, D \,$ as differential 
operator in $\, \d_u \,$.  It cannot be an operator of order zero, 
because $\, B \,$ is a {\it maximal} commutative subalgebra of $\, \D \,$. 
Thus $\, D \,$ has positive order in 
$\, \d_u \,$, so we may apply Proposition~\ref{adbar} 
(with $\, m = r = 0 $) to obtain a $\, q \,$ 
in some extension field of $\, \K \,$ 
such that $\, B \subseteq \c[q] \,$.  The Theorem now follows from 
Lemma~\ref{pol}.
\end{proof}

As in the Introduction, we call a choice of $\, x \,$ as in 
Theorem~\ref{x} a {\it framing} of $\, B \,$.  
Clearly, if $\, x \,$ is a framing of $\, B \,$, then we have 
$\, \Frac B = \c(x) \subseteq \K \,$, where (as in the proof above) 
we have set $\, \K := C_Q(B) \,$.  
In fact it is now easy to see
\begin{theorem}
\label{fr}
If $\, x \,$ is a framing of $\, B \,$, then $\, \c(x) = \K \,$.
\end{theorem}

For the proof of this, we choose the framing $\, x \,$ to be the element 
denoted by $\, u \,$ in the proof of Theorem~\ref{x}, and choose 
$\, \d_x \,$ as in the proof of that Theorem, 
that is, such that $\, [\d_x,x] = 1 \,$ and $\, N_Q(B) = \K[\d_x] \,$. 
We shall think of the elements of $\, \D \,$ as ``operators'' 
(with coefficients in $\, \K \,$). 
\begin{lemma}
\label{beta}
Let $\, L \in \D \,$ (considered as an element of $\, \K[\d_x] \,$) 
have leading coefficient $\, \beta \in \K \,$. Then $\, \beta \in \c(x) \,$.
\end{lemma}
\begin{proof}
By induction on the order of $\, L \,$.  If  $\, L \,$ has  
order zero, that is, $\, L \in \K \,$, then $\, L \,$ commutes with the 
elements of $\, B \,$.  Since $\, B \,$ is a {\it maximal} 
commutative subalgebra of $\, \D \,$, that shows that $\, L \in B \,$, so  
in this case the Lemma just claims that 
$\, B \subset \c(x) \,$, which is certainly true.  Now 
suppose inductively that the assertion is true for 
operators of order $\, n-1 \,$, and let $\, L \in \D \,$ have leading term 
$\, \beta \d_x^{n} \,$.  Fix any $\, b \in B \setminus \c \,$; then 
$\, [L,b] \,$ belongs to $\, \D \,$ and has leading term 
$\, n \beta (\d{b}/\d{x}) \d_x^{n-1} \,$.  
Hence $\, n \beta (\d{b}/\d{x}) \in \c(x) \,$, 
so $\, \beta \in \c(x) \,$.
\end{proof}
\begin{proof}[Proof of Theorem~\ref{fr}]
Let $\, f \in \K \,$.  By the assumption \eqref{ass1}, we have 
$\, f L_1 = L_2 \,$ for some $\, L_i \in \D \,$.  So 
$\, f \beta_1 = \beta_2 \,$, where $\, \beta_i \,$ is the leading 
coefficient of $\, L_i \,$, and the result follows from Lemma~\ref{beta}.
\end{proof}

If $\, x \,$ is \ framing of $\, B \,$ and $\, [\d_x,x] = 1 \,$, 
we shall call the pair $\, (x,\d_x) \,$ a {\it fat framing} of 
$\, B \,$.  Thus we have shown that a fat framing always exists, 
and we have inclusions 
\begin{equation}
\label{incs}
B \subset \D \subset N_Q(B) = \c(x)[\d_x] \subset Q \,.
\end{equation}
The following theorem is a much stronger version of Lemma~\ref{beta}.
\begin{theorem}
\label{y}
Let $\, (x,\d_x) \,$ be a fat framing of $\, B \,$, 
and let $\, L \in \D \,$ be written as an 
element of $\, \c(x)[\d_x] \,$, using the corresponding embedding 
\eqref{incs}. Then the leading coefficient of $\, L \,$ belongs to 
$\, \c[x] \,$.
\end{theorem}

In what follows, for each 
$\, \lambda \in \c \,$ we denote by $\, \mathfrak{v_{\lambda}} \,$ 
the corresponding valuation of $\, \c(x) \,$; that is, if the Laurent 
expansion at $\, \lambda \,$ of a rational function $\, f \,$ has the form
$$
f(x) = \alpha (x-\lambda)^k +\ (\text{higher degree terms})
$$
(with $\, \alpha \not= 0 \,$), then $\, \mathfrak{v_{\lambda}}(f) = k \,$.
Note that $\, \mathfrak{v_{\lambda}}(f') = k-1 \,$, provided 
$\, k \not= 0 \,$.
\begin{lemma}
\label{p}
Let $\, D \in \c(x)[\d_x] \,$ have leading term $\, a(x) \d_x^n \,$, 
where $\, n > 0 \,$, and suppose that $\, D \,$ acts ad-nilpotently on the 
rational function $\, p(x) \,$.  Fix any $\, \lambda \in \c \,$, 
and set
$$
r :=  \mathfrak{v_{\lambda}}(a) \,, \ \ 
s :=  \mathfrak{v_{\lambda}}(p) \,.
$$
Suppose that $\, s \not= 0 \,$.
Then $\, ns = i(n-r) \,$ for some $\, i \in \mathbb{N} \,$. 
\end{lemma} 
\begin{proof}
Let  $\, (\ad D)^i(p) = p_i(x) \d_x^{i(n-1)} + $ (lower order terms), 
so that $\, p_0 = p \,$, $\, p_1 = np'a \,$, and
$$
p_{i+1} = n a p'_i  - i(n-1) a' p_i  \ \ \text{for} \ i \geq 1 \,.
$$
If $\, \mathfrak{v_{\lambda}}(p_i) := q \,$, so that 
$$
p_i = \alpha(x-\lambda)^q + \ldots\ , \quad 
a = \beta(x-\lambda)^r + \ldots\ ,
$$
where $\, \alpha \,$ and $\, \beta \,$ are nonzero and the $\, \ldots \,$ 
denote terms of higher degree in $\, x - \lambda \,$, we find
$$
p_{i+1} = \alpha\beta \{ nq - ir(n-1) \} (x-\lambda)^{q+r-1} + \ldots \,.
$$
So for each $\,i \,$, either 
$\, \mathfrak{v_{\lambda}}(p_{i+1}) = 
\mathfrak{v_{\lambda}}(p_{i}) + r - 1 \,,$
or $\, nq - ir(n-1) = 0 \,$.  Since $\, D \,$ is ad-nilpotent on $\, p \,$,
the latter must occur for some $\, i \,$: let $\, i \,$ 
now denote the {\it first} number 
for which it occurs.  The assumption $\, s \not= 0 \,$ implies that 
$\, \mathfrak{v_{\lambda}}(p_{1}) = r+s-1 \,$ and
$$
q = \mathfrak{v_{\lambda}}(p_{i}) = s + i(r-1)
$$
so $\, n[s + i(r-1)] = ir(n-1) \,$, which simplifies to give the Lemma. 
\end{proof}
\begin{corollary}
\label{ainc}
Let $\, D \in \c(x)[\d_x] \,$ have leading term $\, a(x) \d_x^n \,$, 
where $\, n > 0 \,$.  Suppose $\, D \,$ acts ad-nilpotently on some 
algebra $\, B \subseteq \c[x] \,$ with 
$\, \overline{B} = \c[x] \,$. Then $\, a \in \c[x] \,$. 
\end{corollary}
\begin{proof}
Fix $\, \lambda \in \c \,$.  For any  
$\, s \gg 0 \,$, the algebra $\, B \,$ contains a polynomial 
$\, p \,$ with $ \, \mathfrak{v_{\lambda}}(p) = s \,$ . 
Applying Lemma~\ref{p} to any such $\, p \,$, we find  
that $\, r < n $.  Then applying the lemma with two consecutive values 
of $\, s \,$ and subtracting, we find that $\, n/(n-r) \in \mathbb{N} \,$, 
in particular $\, r \geq 0 \,$.  This shows that $\, a \,$ is regular at 
every point $\, \lambda \in \c \,$, that is, $\, a \,$ is a polynomial.
\end{proof}
Now we can give the 
\begin{proof}[Proof of Theorem~\ref{y}]
As in the proof of Theorem~\ref{x}, we fix an ad-nilpotent element  
$\, D \in \D \setminus B \,$; let its leading term be 
$\, a(x) \d_x^n \,$, where $\, n > 0 \,$. Let $\, L \in \D \,$ 
have leading coefficient $\, \beta(x) \,$.  Then if $\, q, \alpha \,$ 
are as in Proposition~\ref{adbar}, we have $\, \beta \in \c[q, \alpha] \,$. 
We have $\, \alpha^n = a \,$, and by Corollary~\ref{ainc}, $\, a \in \c[x] \,$; 
hence $\, \alpha \,$ is integral over $\, \c[x] \,$.  Also, 
$\, q \,$ is integral over $\, \c[x] \,$ (for if $\, x \,$ has degree 
$\, t \,$ as a polynomial in $\, q \,$, then 
$\, \{ 1, x, \ldots, x^{t-1} \} \,$ generate $\, \c[q] \,$ as  
$\, \c[x]$-module).  Hence every element of  $\, \c[q, \alpha] \,$, in 
particular $\, \beta \,$, is integral over $\, \c[x] \,$. But 
$\, \beta \in \c(x) \,$, hence $\, \beta \in \c[x] \,$.
\end{proof}
\begin{remark}
We have not used the assumption that 
$\, Q \,$ is the Weyl quotient field, so the results of 
this section would still be valid for any of the fields $\, Q \,$ studied 
in Section~\ref{diffop}.  However, this extra generality would be 
illusory, because these fields do not contain any subalgebras 
$\, \D \, $ satisfying \eqref{ass1} and \eqref{ass2} (see \cite{KM}). 
\end{remark}

\section{The dual mad subalgebra}
\label{dual}

The main aim of this section is to show that if a 
certain finiteness condition (\eqref{ass4} below) is 
satisfied, then the mad subalgebra  $\, B \,$ of $\, \D \,$ possesses 
a {\it dual} mad subalgebra  $\, \Check{B} \,$.  If 
$\, \D \,$ is the Weyl algebra $\, \c[x,\d_x] \,$ and 
$\, B = \c[x] \,$, then $\, \Check{B} = \c[\d_x] \,$: in general, the 
relationship between $\, B \,$ and $\, \Check{B} \,$ is similar to this, 
but $\, \Check{B} \,$ is not necessarily isomorphic to $\, B \,$. 

We retain the assumptions  of the preceding 
section; thus $\, \D \,$ is an algebra satisfying 
\eqref{ass1} and \eqref{ass2}, $\, B \,$ is a mad subalgebra of 
$\, \D \,$, and  $\, (x, \d_x) \,$ is a fat framing of  $\, B \,$, so 
that $\, \D \,$ becomes a subalgebra of $\, \c(x)[\d_x] \,$, as in 
\eqref{incs}.  The filtration $\, \{ \D_{\bullet} \} \,$ induced on 
$\, \D \,$ by the usual 
filtration (by order in $\, \d_x $) on $\, \c(x)[\d_x] \,$ 
coincides with that induced by $\, B \,$; in particular, it is 
independent of the choice of fat framing. 
We regard the associated graded algebra
$$
\gr_{\d} \D  :=  \bigoplus_{k \geq 0} \D_k / \D_{k-1}
$$
as embedded in $\, \c(x)[\xi] \,$ via the symbol map (if
$\, L \in \D \,$ has leading term $\, a(x) \d_x^k \,$, its 
symbol is $\, a(x) \xi^k \,$). According to Theorem~\ref{y}, 
we have 
\begin{equation}
\label{ass3}
\gr_{\d} \D \subseteq \c[x, \xi] \,.
\end{equation}  
Following \cite{P}, we now consider the 
{\it $x$-filtration}  on $\, \c(x)[\d_x] \,$ 
(and the filtration it induces on $\, \D \,$). By definition, 
an operator $\, \sum a_i(x) \d_x^i \,$ has $x$-filtration $\, \leq k \,$ 
if $\, \deg_x a_i \leq k \,$ for all $\, i \,$ 
(if $\, f \,$ and $\, g \,$ are polynomials in $\, x \,$, we define 
$\, \deg_x (f/g) := \deg_x f - \deg_x g \,$).  We identify the 
associated graded algebra $\, \gr_x \c(x)[\d_x] \,$ with 
$\, \c[x, x^{-1}, \xi] \,$, and we regard $\, \gr_x \D \,$ as embedded 
in $\, \c[x, x^{-1}, \xi] \,$ via the ``$x$-symbol map'' 
(defined in the obvious way). Theorem~\ref{y} shows that in fact 
$\, \gr_x \D \subseteq \c[x, \xi] \,$, in particular, that the induced 
$x$-filtration on $\, \D \,$ is positive.  We define
\begin{equation}
\label{Bhat}
\check{B} := \{ D \in \D : \deg_x D = 0 \} \ .
\end{equation}
\begin{proposition}
Either $\, \check{B} = \c \,$ or $\, \check{B} \,$ is a mad subalgebra 
of $\, \D \,$.
\end{proposition} 
\begin{proof}
If $\, \check{B} \not= \c \,$, it is easy to check that the 
$x$-filtration is a mad filtration on $\, \D \,$ (see \cite{P}, p.\ 286).
\end{proof}

The following example shows that the undesirable case 
$\, \check{B} = \c \,$ can indeed occur.
\begin{example}
Let $\, \D \subset \c[x, \d_x] \,$ be the subalgebra of the Weyl 
algebra consisting of all operators that can be written as polynomial 
differential operators in the variable $\, w := x^{1/2} \,$.  Then 
$\, \D \,$ contains $\, x \ (=w^2) \,$ and 
$\, x \d_x \ (= \tfrac{1}{2} w \d_w) \,$, so clearly $\, \D \,$ 
satisfies \eqref{ass1}.  Also, $\, \D \,$ contains the mad subalgebra 
$\, B := \c[x] \,$, and the ad-nilpotent element 
$\, \d_w^2 = 4x \d_x^2 + 2 \d_x \,$, hence $\, \D \,$ satisfies 
\eqref{ass2}. But $\, \D \,$ contains no operator (of positive order) 
with constant coefficients, hence $\, \check{B} = \c \,$.
\end{example}

To exclude the possibility that $\, \check{B} = \c \,$, we make 
one more (very strong) assumption about the pair 
$\, (\D,B) \,$, namely
\begin{equation}
\label{ass4}
\gr_{\d} \D \text{\, has finite codimension in \,} \c[x, \xi] \,.
\end{equation}
\begin{proposition}
\label{lead}
Suppose that $\, \D \subseteq \c(x)[\d_x] \,$ satisfies \eqref{ass1}, 
\eqref{ass3} and \eqref{ass4}.  Then the $x$-symbol map defines an 
isomorphism from $\, \check{B} \,$ onto a subalgebra of 
finite codimension in $\, \c[\xi] \,$. In particular,
$\, \check{B} \not= \c \,$. 
\end{proposition}
\begin{proof}
This follows from \cite{LM}, Proposition 2.4, which shows that 
$\gr_{x} \D \,$ and $\gr_{\d} \D \,$ 
have the {\it same finite} codimension in $\, \c[x, \xi] \,$.
Under the $x$-symbol embedding 
$\gr_{x} \D \hookrightarrow \c[x, \xi] \,$ the elements of $\, \c[\xi] \,$ 
come exactly from $\, \check{B} \,$.  Thus if $\, \check{B} \,$ 
had infinite codimension in $\, \c[\xi] \,$,  then $\gr_{x} \D \,$ 
would have infinite codimension in $\, \c[x, \xi] \,$, contradicting 
\cite{LM}.
\end{proof}

We call $\, \check{B} \,$ the {\it dual} mad subalgebra to 
$\, (B,x,\d_x) \,$.  It does not depend on the choice of framing $\, x \,$. 
Indeed, any other framing has the form $\, ax+b \,$ with 
$\, a,b \in \c \,$ and $\, a \neq 0 \,$, so (despite the 
terminology) the $\, x$-filtration on $\, N_Q(B) \,$, and hence on 
$\, \D \,$, does not depend on this choice.  On the other hand, 
$\, \check{B} \,$ 
does depend on the choice of $\, \d_x \,$: a different choice has the 
form $\, \d_x + q \,$ with $\, q \in \c(x) \,$, and if $\, q \,$ has 
positive degree the corresponding  $\, x$-filtration, and hence also 
$\, \check{B} \,$, may change drastically. However, we do have the following.
\begin{lemma}
\label{nochange}
If we change $\, \d_x \,$ to $\, \d_x + q \,$, where $\, q \in \c(x) \,$ 
has \emph{negative} degree in $\, x \,$, then the $\, x$-filtration on 
$\, N_Q(B) \,$, and hence also the dual mad subalgebra $\, \check{B} \,$, 
remain unchanged.
\end{lemma}

Proposition \ref{lead} implies that $\, \Check{B} \,$ contains an operator 
of every sufficiently high order, that is, that $\, \Check{B} \,$ 
is an {\it algebra of rank} $1$ in the sense of Burchnall-Chaundy theory  
(cf.\ \cite{BC}).
By \eqref{ass3}, the leading coefficient of every operator 
in $\, \Check{B} \,$ is constant; however, in the Burchnall-Chaundy theory
it is convenient to consider algebras of differential 
operators that are normalized to have their first {\it two} 
coefficients constant.  We call a fat framing $\, (x,\d_x) \,$ of 
$\, B \,$ {\it good} if the corresponding $\, \Check{B} \,$ has 
this property.
\begin{proposition}
\label{good}
Let $\, (x,\d) \,$ be any fat framing of $\, B \,$.  Then $\, B \,$ 
has a good fat framing $\, (x,\d_x) \,$ with the same dual subalgebra 
$\, \Check{B} \,$.
\end{proposition}
\begin{proof}
Choose any $\, L \in \Check{B} \,$ of positive order and with leading 
coefficient $\, 1 \,$: it has the form
$$
L = \d^n + (c + nq) \d^{n-1} + \text{ (lower order terms)} 
$$
where $\, c \in \c \,$ and $\, \deg_x q < 0 \,$.  Let 
$\, \d_x := \d + q \,$: then by Lemma~\ref{nochange} the fat framings 
$\, (x,\d_x) \,$ and $\, (x,\d) \,$ determine the same $\, \Check{B} \,$, 
and we have   
$$
L = \d_x^n + c \d_x^{n-1} + \, (\text{lower order terms}) \,; 
$$
that is, the first two coefficients of $\, L \,$ are now constant.  
Any operator that commutes with $\, L \,$ also 
has this property, hence all the elements of $\, \Check{B} \,$ 
now have their first two coefficients constant; that is, 
$\, (x,\d_x) \,$ is good.
\end{proof}
\begin{remark}
\label{goodness}
The notion of a ``good'' fat framing introduced above may seem 
a little artificial. To appreciate it better, let us 
reconsider the case where $\, \D \,$ is 
the Weyl algebra $\, A \,$.  By Theorem~\ref{D1}, 
in this case we have $\, x \in A \,$ for any framing $\, x \,$ of 
a mad subalgebra; a fat framing $\, (x,\d_x) \,$ is good 
exactly when $\, \d_x \in A \,$ too. It follows easily 
from Theorem~\ref{D2} that the group $\, \Aut A \,$ acts {\it freely and 
transitively} on the set of triples $\, (B,x,\d_x) \,$, where 
$\, B \,$ is a mad subalgebra of $\, A \,$ and 
$\, (x,\d_x) \,$ is a good fat framing of $\, B \,$.  
We shall see 
in Section~\ref{10} that the same is true for any of our algebras 
$\, \D \,$, except that $\, \Aut \D \,$ has to be replaced by the 
larger group $\, \Pic \D \,$ (in the case of the Weyl algebra these two 
groups coincide).
\end{remark}

\section{The adelic Grassmannian}

In this section we summarize various facts about the adelic 
Grassmannian $\, \Gr $ which we need to prove our main results.  
We make no attempt to indicate proofs, except for Theorem~\ref{bingr}, 
which we have not been able to find stated explicitly in the literature.

\subsection{The Grassmannian}
\label{grdef}
We recall the definition of 
$\, \Gr$.  For each $\, \lambda \in \c \,$, 
we choose a {\it $\lambda$-primary} subspace of 
$\, \c[z] \, $, that is, a linear subspace  $\, V_{\lambda} \,$
such that 
$$
(z - \lambda)^{N} \c[z] \subseteq V_{\lambda} \ \text{for some} \ N \, .
$$
We suppose that $\, V_{\lambda} = \c[z] \, $ for all but finitely 
many $\, \lambda \,$.  Let $\, V = \bigcap_{\lambda} V_{\lambda} \,$ 
(such a space $\, V \,$ is called {\it primary decomposable}) and, 
finally, let
$$
W = \prod_{\lambda} (z - \lambda)^{- k_{\lambda}} \, V \subset 
\c(z) \, ,
$$
where $\, k_{\lambda} \,$ is the codimension of $\, V_{\lambda} \,$ 
in $\, \c[z] \,$. By definition, $\, \Gr$ consists of all 
$\, W \subset \c(z) \,$ obtained in this way.  For each 
$\, W \in \Gr$ we set
\begin{equation}
\label{aw}
A_W := \{ f \in \c[z] : fW \subseteq W \} \,;
\end{equation}
then the inclusion $\, A_W \subseteq \c[z] \,$ corresponds to a 
framed curve $\, \pi : \A^1 \to X \,$ and the $\, A_W$-module $\, W \,$ 
corresponds to a rank 1 torsion-free coherent sheaf $\, \L \,$ over 
$\, X \,$. Indeed, in this way the points of $\, \Gr $ correspond 
bijectively to isomorphism classes of such triples $\, (\pi, X, \L) \,$. 
For more details see \cite{W1}.

\subsection{The Baker function and the Burchnall-Chaundy theory}
Associated to each $\, W \in \Gr$ is its {\it Baker function} 
$\, \psi_W \,$ (see \cite{SW} or \cite{W1}). It has the form
\begin{equation}
\label{bakad}
\psi_W(x,z) = e^{xz} \{ 1 + \sum_i f_i(x) g_i(z) \} \ ,
\end{equation}
where the $\, f_i, g_i \,$ are rational functions that vanish at infinity. 
For each $\, f \in A_W \,$ there is a unique differential operator 
$\, L_f \in \c(x)[\d_x] \,$ such that 
\begin{equation}
\label{Lf}
L_f \psi_W(x,z) = f(z) \psi_W(x,z) \ ;
\end{equation}
the map $\, f \mapsto L_f \,$ defines an isomorphism  from $\, A_W \,$ 
to a maximal commutative rank 1 subalgebra $\, \mathcal{A}_W \,$ of 
$\, \c(x)[\d_x] \,$. Clearly, the operators $\, L_f \,$ are normalized 
to have their first two coefficients constant. 

We shall need (briefly) the larger Grassmannian $\, \Grat $ of 
\cite{W1}: it is similar to $\, \Gr $, except that the normalization 
map $\, \pi : \A^1 \to X \,$ is not required to be bijective.  The 
Baker function now does not necessarily have the form \eqref{bakad}, 
and the operators $\, L_f \,$ may not have rational coefficients; however, 
we can expand $\, \psi_W  \,$ in a series
\begin{equation}
\label{exp}
\psi_W(x,z) = e^{xz} \{ 1 + \sum_1^{\infty} a_i(x) z^{-i} \}
\end{equation}
in which the coefficients $\, a_i \,$ are rational functions of 
$\, x \,$ and some exponentials $\, e^{\lambda x} \,$ (the 
numbers $\, \lambda \,$ occurring are the inverse images under 
$\, \pi \,$ of the multiple points of $\, X \,$). 
Every normalized rank 1 algebra of differential operators 
$\, \mathcal{A} \,$ with $\, \Spec \mathcal{A} \,$ rational can be 
obtained from a point of $\, \Grat $ in the way explained above for 
$\, \Gr $. The following Theorem is almost proved in \cite{W1}.
\begin{theorem}
\label{bingr}
Let $\, \B \subseteq \c(x)[\d_x] \,$ be any rank {\rm 1} commutative 
algebra of differential operators with first two coefficients 
constant, and\,\footnote{This last assumption is almost certainly superfluous, 
but we do not know a reference.}
such that the curve $\, \Spec \B \,$ is 
rational. 
Then there is a unique $\, W \in \Gr $ such that 
$\, \B \subseteq \mathcal{A}_W \,$.
\end{theorem}
\begin{proof}
Let $\, \mathcal{A} \,$ be the maximal commutative algebra of differential 
operators containing $\, \B \,$.  Then $\, \Spec \mathcal{A} \,$ 
is still a rational curve (with normalization $\, \A^1 $), so it is known 
that $\, \mathcal{A} = \mathcal{A}_W \,$ for some $\, W \in \Grat \,$.  
The assertion that $\, W \,$ can be chosen to be in $\, \Gr $ is 
equivalent to saying that $\, \Spec \mathcal{A} \,$ is a framed curve. 
According to \cite{W1}, that in turn is equivalent to the fact that 
if the Baker function of $\, W \,$ is expanded in the form \eqref{exp},
then all the $\, a_i \,$ are rational functions of $\, x\,$.  But 
if $\, \mathcal{A}_W \,$ contains an operator of positive 
order with rational coefficients, then this must be the case. 
To see that, let 
$$
L = \d_x^n + c \d_x^{n-1} + \sum_0^{n-2} u_i(x) \d_x^i  \quad 
(n>0, \ c \in \c)
$$
be such an operator; then we have an equation
\begin{equation}
L \psi_W(x,z) = (z^n + c z^{n-1} + \ldots) \psi_W(x,z) \ .
\end{equation}
Substituting in the expansion \eqref{exp} of $\, \psi_W \,$ and 
equating coefficients of powers of $\, z \,$, we get a recursion 
relation of the form
$$
a_r' = \{ \text{some polynomial in derivatives of the $\, u_i \,$ 
and the $\, a_j \,$ with $\, j < r \,$} \} \ .
$$
Now suppose inductively that the $\, a_j \,$ are rational for 
$\, j < r \,$. The recursion relation then shows that  $\, a_r \,$ 
is the sum of a rational function and (possibly) some logarithmic 
terms $\, \lambda \log (x - \mu) \,$.  But $\, a_r \,$ is a 
meromorphic function, hence the logarithmic terms must be absent.  
Thus all $\, a_i \,$ are rational, as claimed.
\end{proof}
\subsection{The algebras $\, \D(W) \,$}
\label{DW}

For each $\, W \in \Gr $, we define
\begin{equation}
\D(W) := \{ D \in \c(z)[\d_z] : D.W \subseteq W \} \ .
\end{equation}
If $\, W \,$ corresponds to the triple $\, (\pi, X, \L) \,$, then 
we can interpret $\, \D(W) \,$ as the ring 
$\, \D_{\L}(X) \,$ of differential operators 
on sections of the sheaf $\, \L \,$ (embedded in $\, \c(z)[\d_z] \,$ 
via the ``framing'' $\, \pi $).
It is fairly well-known (cf.\ \cite{SS}, or see the Appendix) 
that the algebras $\, \D(W) \,$ 
satisfy all the assumptions we have made about $\, \D \,$ in the 
previous sections: we shall prove later (see Corollary~\ref{wow})
that in fact the $\, \D(W) \,$ 
are (up to isomorphism) the {\it only} algebras that satisfy these 
assumptions.

The paper \cite{CH} provides further information about the algebras 
$\, \D(W) \,$: let us list the results that we need from 
that paper.  If $\, V \,$ and $\, W \,$ are linear subspaces of 
$\,\c(z) \,$, we set 
\begin{equation}
\D(V,W) := \{ D \in \c(z)[\d_z] : D.V \subseteq W \} \,.
\end{equation}
If $\, V, \,W \in \Gr $, then clearly  $\, \D(V,W) \,$ is a 
$\, \D(W)$-$\D(V)$-bimodule (the actions being given by composition). 
\begin{theorem}[\cite{CH}]
\label{CH1}
Each isomorphism class of right ideals in the Weyl algebra 
$\, A \equiv \c[z,\d_z] \,$ has a unique representative of the form 
$\, \D(\c[z], W) $ with $\, W \in \Gr $.
\end{theorem}

Generalizing slightly the definition in Subsection~\ref{grdef}, let 
us call a linear subspace $\, V \subset \c(z) \,$ primary decomposable if 
$\, V = fW \,$ for some $\, f \in \c(z) \,, \ W \in \Gr $.  
\begin{theorem}[\cite{CH}]
\label{CH2}
A subspace $\, V \subset \c(z) \,$ is primary decomposable if and only if 
$\, \D(\c[z],V).{\c[z]} =  V \,$.
\end{theorem}

\begin{theorem}[\cite{CH}]
\label{CH3}
For each $\, W \in \Gr $, the algebra $\, \D(W) \,$ can be identified 
with the endomorphism ring of the corresponding $\, A$-module 
$\, \D(\c[z],W) \,$.
\end{theorem}

Since the Weyl algebra $\, A \,$ is hereditary and simple, 
every ideal in it is 
a progenerator; so Theorem~\ref{CH3} implies that all the algebras 
$\, \D(W) \,$ are Morita equivalent to $\, A \,$; in particular, 
all the $\, \D(W) \,$ are {\it simple}.  Furthermore, all the 
bimodules $\, \D(V,W) \,$ are invertible, and for any 
$\, U,V,W \in \Gr $, we have
\begin{equation}
\label{uvw}
\D(V,U) \D(W,V) = \D(W,U) \,.
\end{equation}

\subsection{The action of $\, \Gamma \,$}
\label{gammact}

In the theory of integrable systems, a key role is played by the 
action on $\, \Gr $ of the group $\, \Gamma \,$ from the  
Introduction.  Recall that for each polynomial $\, p(z) \,$ we have 
the element $\, \gamma_p \in \Gamma \,$ defined by \eqref{gamma}: 
it acts on the Weyl algebra, or more generally on the algebra 
$\, \c(z)[\d_z] \,$ as formal conjugation by $\, e^{p(z)} \,$.  
Roughly speaking, the action of $\, \Gamma \,$ on $\, \Gr $ is 
given by scalar multiplications; that is, 
we define $\, \gamma_p W := e^{p(z)} W \,$.  Of course, since 
$\, e^{p(z)} \,$ is not a rational function, this does not 
immediately make sense: to interpret it correctly
we have temporarily to replace $\, W \,$ by 
a suitable completion (see, for example \cite{BW2}, Section 2). 
This difficulty need not concern us here, because 
we are interested mainly in the induced action of $\, \Gamma \,$ on 
the spaces $\, \D(V,W) \,$, which makes sense without 
any completions.  Namely, we have 
$$
\D(\gamma_p^{-1} V, \gamma_p^{-1} W) = e^{-p(z)} \D(V,W) e^{p(z)} \ , 
$$
so that $\, \gamma_p \,$ induces a bijective map 
$$
\gamma_p : \D(\gamma_p^{-1} V, \gamma_p^{-1} W) \to \D(V,W) 
$$
defined by $\, \gamma_p (D) := e^{p(z)} D e^{-p(z)} \,$.  In particular, 
taking $\, V = W \,$, we have isomorphisms of algebras
$$
\gamma : \D(\gamma^{-1} W) \to \D(W) 
$$
for each $\, \gamma \in \Gamma, \  W \in \Gr $.  We refer to 
\cite{BW2} for a more thorough discussion of these points.

\subsection{The bispectral involution}
\label{bi}
The {\it bispectral involution} 
$\, W \mapsto b(W) \,$ on $\, \Gr $ can be
characterized by the formula 
\begin{equation}
\label{b}
\psi_{b(W)}(x,z) = \psi_{W}(z,x) \ .
\end{equation}
Generalizing \eqref{Lf}, one can show (see \cite{BW2}) 
that for each $\, D \in \D(W) \,$ 
there is a unique differential operator $\, \Theta \,$ in the 
variable $\, x \,$ such that
\begin{equation}
\label{Theta}
D(z).\psi_W(x,z) = \Theta(x).\psi_W(x,z) \ .
\end{equation}
The map $\, D \mapsto \Theta \,$ defines\footnote{After 
restoring the notation $\, z \,$ for $\, x \,$; this kind of confusion will 
recur several times in what follows.} 
an anti-isomorphism 
from $\, \D(W) \,$ to 
$\, \D(b(W)) \,$.  To write it more explicitly, we introduce 
the {\it formal integral operator} (in $\, x $) $\, K_W \,$ 
with the property that (formally) $\, \psi_W = K_W . e^{xz} \,$.  
If $\, \psi_W \,$ is given by \eqref{bakad}, then we have
\begin{equation}
\label{K}
K_W = 1 + \sum_i f_i(x) g_i(\d_x) \ ;
\end{equation}
note that $\, K_W \,$ belongs to the Weyl quotient field $\, Q \,$.  
If we denote by $\, b \,$ also the anti-automorphism of $\, Q \,$ 
defined by $\, b(x) = \d_x , \  b(\d_x) = x \,$, then the formula 
\eqref{b} takes the form
$$
K_{b(W)} = b(K_W) \, ,
$$
and \eqref{Theta} says that the 
anti-isomorphism $\, \beta : \D(W) \to \D(b(W)) \,$
defined above is given by the formula
\begin{equation}
\label{bD}
\beta(D) = K_W b(D) K_W^{-1}\,.
\end{equation}

The connection of the bispectral involution with the construction 
in Section~\ref{dual} is as follows.

\begin{proposition}
The algebra $\, \beta^{-1}(A_{b(W)}) \equiv \mathcal{A}_{b(W)} \,$ 
is the mad subalgebra of $\, \D(W) \,$ dual to $\, A_W \,$.
\end{proposition}
\begin{proof}
This follows at once from \eqref{bD} and the fact  
that $\, K_W - 1 \,$ has negative $x$-filtration.
\end{proof}

\section{Proof of Theorems~\ref{T1} and \ref{T2}}
\label{P1}

We now come back to the situation of Section~\ref{dual}: thus
we have the mad subalgebra $\, B \subset\D \,$ together with 
a good fat framing $\, (x,\d_x) \,$ of $\,B \,$. The dual 
subalgebra $\, \Check{B} \,$ then satisfies 
the conditions of Theorem~\ref{bingr},
so it determines a point of the adelic Grassmannian. We denote this point 
by $\, b(W) \,$ (where $\, b \,$ is the 
bispectral involution on $\, \Gr $) so that  
$\, \Check{B} \subseteq \mathcal{A}_{b(W)} \,$.  As in 
Subsection~\ref{bi}, we allow ourselves the imprecision of using 
$\, x \,$ to denote the variable in the definition of $\, \Gr$, so that 
$\, W \,$ is a subspace of $\, \c(x) \,$, and both $\, \D \,$ and 
$\, \D(W) \,$ are subalgebras of $\, \c(x)[\d_x] \,$.  
With that understanding, the main 
result of this section is as follows.
\begin{theorem}
\label{main}
With $\, W $ defined as above, we have $\, \D = \D(W) \,$.
\end{theorem}
\begin{proof}
For each $\, L \in \Check{B} \,$ 
we have an equation of the form 
$\, L \psi_{b(W)} = f(z) \psi_{b(W)} \,$, or equivalently, 
$\, LK_{b(W)} = K_{b(W)} f(\d_x) \,$ 
($\, f\,$ is a polynomial). 
Thus $\, K_{b(W)}^{-1} \Check{B} K_{b(W)} \,$ is a subalgebra of 
$\, \c[\d_x] \,$. 
Since $\, \Check{B} \,$ acts ad-nilpotently on $\, \D \,$, 
the algebra 
$\, K_{b(W)}^{-1} \Check{B} K_{b(W)} \,$ acts ad-nilpotently on 
$\, K_{b(W)}^{-1} \D K_{b(W)} \,$, hence
$$
K_{b(W)}^{-1} \D K_{b(W)} \subseteq N_Q(\d_x) = \c(\d_x)[x] \ .
$$
Applying the anti-involution $\, b \,$, we deduce that 
$$
K_{W} b(\D) K_{W}^{-1} \subseteq \c(x)[\d_x] \ ;
$$
by \cite{BW2}, Proposition 8.2, that is equivalent to 
$$
\D \subseteq \D(W) \ .
$$
To see that we have equality here, we use the following 
lemma of Levasseur and Stafford (see \cite{LS}): 
let $\, R \subseteq S \,$ 
be Noetherian domains such that 
(i)  $\, R \,$ and $\, S \,$ have the same quotient field; 
(ii) one of  $\, R \,$ and $\, S \,$ is simple;
(iii)  $\,S \,$ is finitely generated as an $\, R$-module 
(both left and right). 
Then $\, R = S \,$.
Let us check that the hypotheses of the lemma 
are satisfied for $\, \D \subseteq \D(W) \,$.  
Certainly, these are both domains with quotient field $\, Q \,$, 
and $\, \D(W) \,$ is simple.  Because  the finiteness condition 
\eqref{ass4} is satisfied,  $\, \c[x, \xi] \,$ is a 
finitely generated module over 
$\, \gr_{\d} \D \,$ (or $\, \gr_{\d} \D(W) \,$), so by  
\cite{AM}, Proposition 7.8, these are finitely generated $\c$-algebras, 
hence Noetherian rings.  It follows that 
$\, \D \,$ and $\, \D(W) \,$ are also (both left and right) Noetherian.  
Finally, to see the property (iii), note that 
we have 
$$
\gr_{\d} \D \subseteq \gr_{\d} \D(W) \subseteq \c[x, \xi] \ ,
$$
and $\, \gr_{\d} \D \,$ has finite codimension in $\, \c[x, \xi] \,$, 
hence ({\it a fortiori}) in $\,\gr_{\d} \D(W) \,$.  Thus 
$\,\gr_{\d} \D(W) \,$ is a finitely generated module over 
$\, \gr_{\d} \D \,$, so property (iii) follows.
\end{proof}

If we now combine Theorem~\ref{main} with the main result of \cite{LM}, 
we obtain the following.
\begin{corollary}
\label{wow}
Let $\, \D \,$ be a subalgebra of the Weyl quotient field $\, Q \,$ 
satisfying \eqref{ass1} and \eqref{ass2}.  Then either all mad 
subalgebras $\, B \,$ of  $\, \D \,$ satisfy 
the finiteness condition \eqref{ass4}, or else none of them does. 
Furthermore, in the former case $\, \D \,$ is isomorphic to $\, \D(X) \,$ 
for some framed curve $\, X \,$. 
\end{corollary}
\begin{proof}
Suppose that $\, \D \,$ possesses one mad subalgebra satisfying 
\eqref{ass4}.  According to Theorem~\ref{main}, this implies 
that $\, \D \,$ is isomorphic to $\, \D(W) \,$ for some 
$\, W \in \Gr$, and hence to $\, \D(X) \,$ for some framed curve 
$\, X \,$ (see \cite{BW1} or \cite{BW3}).  
It is well known (see \cite{SS}) that 
the pair $\, (\D(X), \O(X)) \,$ satisfies \eqref{ass4}; 
the main result of 
\cite{LM} states that for the algebras $\, \D(X) \,$ the codimension 
of $\, \gr_{\d} \D \,$ in $\, \c[x, \xi] \,$ is independent of the 
choice of mad subalgebra $\, B \,$; in particular, it is always 
finite, as claimed.  That completes the proof.
\end{proof}

It is now easy to give the 
\begin{proof}[Proof of Theorems~\ref{T1} and \ref{T2}]
Let $\, \D(X) \subset \c(z)[\d_z] \,$ be the ring of differential 
operators on a framed curve, and let $\, B \,$ be a mad 
subalgebra of $\, \D(X) \,$. 
By the above, we may choose a 
good fat framing of $\, B \,$; then Theorem~\ref{main} gives us 
a point $\, W \,$ of $\, \Gr $ and an isomorphism 
$\, \varphi: \D(W) \to \D(X) \,$ taking $\, A_W \,$ onto $\, B \,$.  
Thus $\, Y := \Spec B \simeq \Spec A_W \,$ is a framed curve, as 
claimed in Theorem~\ref{T1}.  And Theorem~\ref{T2} follows at once, because 
$\, \D(W) \simeq \D_{\L}(Y) \,$ 
for some sheaf $\, \L \,$ over $\, Y \,$, (see Subsection~\ref{DW}).
\end{proof}

\section{An alternative proof of Theorem~\ref{main}}
\label{P2}
As we have just seen, our Theorems~\ref{T1} and \ref{T2}  
are proved by combining the main result of \cite{LM} 
with Theorem~\ref{main}.   
The proof of Theorem~\ref{main} given in the preceding 
section depends heavily on machinery 
inspired by the theory of integrable systems.  In the present 
section we want to give a proof that makes the minimum possible 
use of this machinery; namely, we shall use from it only the following 
consequence of Theorem~\ref{bingr}.
\begin{proposition}
\label{psi0}
Let $\, \B \subseteq \c(x)[\d_x] \,$ be any rank {\rm 1} commutative 
algebra of differential operators with first two coefficients 
constant, and such that the curve $\, \Spec \B \,$ is rational.
Then there is a rational function $\, \psi_0(x) \,$ whose annihilator in 
$\, \B \,$ is a maximal ideal of $\, \B \,$.
\end{proposition}
\begin{proof}
Let $\, \psi(x,z) \,$ be a joint eigenfunction for $\, \B \,$ of 
the form \eqref{bakad}, so that for each $\, L \in \B \,$ 
we have an equation
$\, L \psi(x,z) = f_L(z) \psi(x,z) \,$.
Suppose first that $\, \psi \,$ is regular at $\, z = 0 \,$, and set 
$\, \psi_0 := \psi(x,0) \,$.  Then $\, \psi_0 \in \c(x) \,$, and 
$\, L \psi_0 = f_L(0) \psi_0 \,$ for all  $\, L \in \B \,$, 
so the annihilator of $\, \psi_0 \,$ is the kernel of the character 
$\, L \mapsto f_L(0) \,$ of $\, \B \,$. If $\, \psi \,$ 
has a pole of order $k$ at $\, z=0 \,$ we replace it by 
$\, z^k \psi \,$ and argue as above.
\end{proof}

Returning to our algebra $\, \D \,$ with its pair of mad subalgebras 
$\, (B, \Check{B}) \,$, we may apply Proposition~\ref{psi0} to the rank 1 
algebra $\,\Check{B} \,$: let $\, V := \D.\psi_0 \,$ be the cyclic 
sub-$\D$-module of 
$\, \c(x) \,$ generated by the corresponding function $\, \psi_0 \,$. 
We aim to show that $\, V \,$ coincides with the space 
$\, W \in \Gr$ of the preceding section. In contrast to what we had 
there, it is clear from the definition of $\, V \,$ that 
$\, \D \subseteq \D(V) \,$; however, it is not clear that  
$\, V \in \Gr $.  The crucial step towards proving that is the following.
\begin{lemma}
\label{Vfinite}
$ V $ is finite over $\, B \,$.
\end{lemma}
\begin{proof}
Let $\, I \subset \D \,$ be the annihilator 
of $\, \psi_0 \,$ in $\, \D \,$, and let 
$\, \mathfrak{m} = I \cap \Check{B} \,$: according to 
Proposition~\ref{psi0}, $\, \mathfrak{m} \,$ is a  
maximal ideal in $\, \Check{B} \,$.
Clearly, $\, I \,$ contains the extension $\, \D \mathfrak{m} \,$ of 
$\, \mathfrak{m} \,$ to $\, \D \,$ (in fact $\, I = \D \mathfrak{m} \,$, 
but we do not need to prove that here).  Thus 
$\, V \simeq \D/I \,$ is a quotient module of $\, \D/{\D \mathfrak{m}} \,$, 
so it is enough to prove that $\, \D/{\D \mathfrak{m}} \,$ is finite over 
$\, B \,$. We regard $\, \D/{\D \mathfrak{m}} \,$ as a filtered 
$\D$-module (via the $x$-filtration): the associated graded module can 
then be identified\footnote{To simplify the notation, we 
do not distinguish between $\, \mathfrak{m} \subset \Check{B} \,$ 
and its isomorphic image in $\, \gr_x \Check{B} \,$.} 
with $\, \gr_x \D/({\gr_x \D) \mathfrak{m}} \,$.  Thus it is enough if we 
prove that this is finite over $\, \gr_x B \,$.  Choose 
$\, p(\xi) \in \gr_x \Check{B} \,$ so that 
$\, p(\xi) \c[x, \xi] \subseteq \gr_x \D \,$ (that is possible, 
because $\, \c[x, \xi] / \gr_x \D \,$ is a finite-dimensional 
$\, \gr_x \Check{B}$-module, so its annihilator is a nonzero ideal in 
$\, \gr_x \Check{B} \,$).  Let 
$\, \mathfrak{n} := \c[\xi] p(\xi) \mathfrak{m} \,$ (thus 
$\, \mathfrak{n} \,$ 
is a nonzero ideal in $\, \c[\xi] \,$).  We have 
\begin{equation}
\label{list}
\c[x, \xi] \mathfrak{n} \equiv \c[x, \xi] p(\xi) \mathfrak{m} \,\subseteq 
\, (\gr_x \D) \mathfrak{m} \,\subseteq \, \gr_x \D \,\subseteq \, \c[x, \xi] \ .
\end{equation}
Now, $\, M := \c[x, \xi] / \c[x, \xi] \mathfrak{n} \,$ is a finite 
$\c[x]$-module (in fact it is free of rank equal to the codimension 
of $\, \mathfrak{n} \,$ in $\, \c[\xi] \,$), and $\, \c[x] \,$ 
is a finite $\, \gr_x B$-module (because $\, \gr_x B \,$ has finite 
codimension in $\, \c[x] \,$).  Thus $\, M \,$ is a finite 
$\, \gr_x B$-module, and hence Noetherian (because $\, \gr_x B \,$ 
is Noetherian).  Thus the subquotient (see \eqref{list})  
$\, \gr_x \D / (\gr_x \D) \mathfrak{m} \,$ of $\, M \,$ is 
again a Noetherian $\, \gr_x B$-module.
\end{proof}

For the rest of this section $\, V \,$ could be any sub-$\D$-module 
of $\, \c(x) \,$ that is finite\footnote{Readers who wish to avoid using 
the Burchnall-Chaundy theory have only to prove the existence of 
such a $\, V \,$.} over $\, B \,$: however, we shall see at the end of this 
section that $\, V \,$ is in fact uniquely determined by these properties.
Being finite over $\, B \,$ means that $\, V \,$ is the space of sections 
of a rank 1 torsion-free sheaf over the curve 
$\, \Spec B \,$. The next Proposition 
is thus (part of) Proposition 7.1 in \cite{SW}, but we shall give a 
self-contained proof. 
\begin{proposition}
\label{grat}
Let $\, V \subset \c(x) \,$ be as above.  Then there are nonzero 
polynomials $\, p, q \in \c[x] \,$ such that
$\, pV \subseteq \c[x] \text{\, and \,} q \c[x] \subseteq V \,$.
\end{proposition}
\begin{proof}
Let $\, \{ v_1, \ldots, v_k \} \,$ generate $\, V \,$ as a $B$-module, 
and let $\, v_i = f_i/g_i \,$, where $\, f_i \,$ and $\, g_i \,$ are 
polynomials.  Then if $\, p \,$ is the product (or least common multiple) 
of the $\, g_i \,$, clearly $\, pv_i \in \c[x] \text{\\, for all \,} i \,$, 
hence $\, pV \subseteq \c[x] \,$.  Now let  
$\, \{ f_1, \ldots, f_r \} \,$ generate $\, \c[x] \,$ as a $B$-module, 
and let $\, f_i = b_i/c_i \,$, where $\, b_i, c_i \in B \,$.  If 
$\, a \,$ is the product (or least common multiple) 
of the $\, c_i \,$, then $\, af_i \in B \text{\\, for all \,} i \,$, 
hence $\, a \c[x] \subseteq B \,$.  On the other hand, let $\, v \,$ 
be any nonzero element of $\, V \,$, and let $\, v = b/c \,$, where 
$\, b, c \in B \,$.  Then $\, cv = b \in V \cap B \,$; since $\, V \,$ 
is a $B$-module, it follows that $\, Bb = bB \subseteq V \,$.  Thus 
if $\, q := ba \,$ then we have
$\, q \c[x] = ba \c[x] \subseteq bB \subseteq V $.
\end{proof}

Now set
\begin{equation}
\label{P}
P := \D(\c[x], V) \equiv \{ D \in \c(x)[\d_x] : D.\c[x] \subseteq V \} \ .
\end{equation}
Clearly, $\, P \,$ is a right sub-$A$-module of $\, \c(x)[\d_x] \,$.
\begin{corollary}
\label{ideal}
$\, P \,$ is a (fractional) right ideal of $\, A \,$ having nonzero 
intersection with $\, \c(x) \,$.
\end{corollary}  
\begin{proof}
Let $\, p \,$ be as in Proposition~\ref{grat}.  Then 
$\, pP.\c[x] \subseteq pV \subseteq \c[x] \, $, 
hence $\, pP \subseteq A \,$, so $\, pP \,$ is an (integral) 
right ideal in $\, A \,$.  And if $\, q \,$ is as in 
Proposition~\ref{grat}, then by definition $\, q \in P \,$; thus 
$\, P \,$ has nonzero intersection even with $\, \c[x] \,$.
\end{proof}

Finally, we set 
\begin{equation}
\label{E}
E := \{ D \in \c(x)[\d_x] : DP \subseteq P \} \ .
\end{equation}
\begin{lemma}
\label{end}
We may identify $\, E \,$ with the endomorphism ring $\, \End_A P \,$.
\end{lemma}
\begin{proof}
Every $A$-module endomorphism of $\, P \,$ extends uniquely to a 
$Q$-linear endomorphism of the $1$-dimensional (right) vector space 
$\, P \otimes_A Q \simeq Q \,$: it follows that we may (as is usual) 
identify $\, \End_A P \,$ with the algebra
$$
E' := \{ D \in Q: DP \subseteq P \} \ .
$$
By Corollary~\ref{ideal}, we may choose a nonzero element 
$\, q \in P \cap \c(x) \,$; then if 
$\, D \in E' \,$, we have $\, Dq \in P \subseteq \c(x)[\d_x] \,$, 
hence $\, D \in \c(x)[\d_x] \,$. Thus $\, E' = E \,$.
\end{proof}
\begin{lemma}
\label{grE}
$ \gr_{\d} E \subseteq \c[x, \xi] \,$.
\end{lemma}
\begin{proof}
As in the proof of Corollary~\ref{ideal}, we have $\, pP \subseteq A \,$ 
for suitable $\, p \in \c[x] \,$; hence 
$\, p \gr_{\d} P \subseteq \gr_{\d} A = \c[x, \xi] \,$.  
Thus $\, \gr_{\d} P \,$ is a fractional ideal of $\, \c[x, \xi] \,$, and 
hence is a finitely generated  $\, \c[x, \xi]$-module (since this ring is 
Noetherian).  Let $\, \{ p_1, \ldots, p_m \} \,$ generate 
$\,\gr_{\d} P \,$ as a 
$\, \c[x, \xi]$-module, and let $\, d \in \gr_{\d} E \,$.  It follows 
from the definition \eqref{E} of $\, E \,$ that 
$\, d \gr_{\d} P \subseteq \gr_{\d} P \,$, so we have equations of the form
$\, dp_i = \sum_1^m f_{ij} p_j \,$ for some $\, f_{ij} \in \c[x, \xi] \,$. 
Multiplying on the left by the adjoint of the matrix 
$\, (d \delta_{ij} - f_{ij}) \,$, we find that 
$\, \det(d \delta_{ij} - f_{ij})$ annihilates the $\, p_i \,$, 
hence it is zero.  That shows that $\, d \,$ is integral over 
$\, \c[x, \xi] \,$; but this ring is integrally closed, hence 
$\, d \in \c[x, \xi] \,$, as claimed.
\end{proof}

The next result completes our alternative proof of Theorem~\ref{main}.
\begin{theorem}
\label{main1}
We have $\, \D = \D(V) = E \,$.  Moreover, $\, V \in \Gr $.     
\end{theorem}
\begin{proof}
Clearly, we have inclusions of algebras
\begin{equation}
\label{inc}
\D \subseteq \D(V) \subseteq E \ .
\end{equation}
To see that we have equality here we use the Levasseur-Stafford 
lemma, just as in the proof of Theorem~\ref{main}.  
Let us check that the conditions of that lemma are satisfied 
by the pair of algebras $\, \D \subseteq E \,$.
First, since $\, \D \,$ satisfies \eqref{ass1} and $\, E \subset Q \,$, 
it is obvious that $\, \D \,$ and $\, E \,$ both have quotient 
field $\, Q \,$.  Next, 
the algebra $\, A \,$ is hereditary and simple, and $\, P \,$ is an 
ideal of $\, A \,$ (see corollary~\ref{ideal}):
hence $\, P \,$ is a progenerator, so 
by Lemma~\ref{end} $\, E \,$ is Morita equivalent 
to $\, A \,$. It follows that $\, E \,$ is simple and Noetherian.  Finally, 
using \eqref{inc} and Lemma~\ref{grE}, we get inclusions of 
commutative algebras 
\begin{equation}
\label{grinc}
\gr_{\d} \D \subseteq \gr_{\d} \D(V) \subseteq 
\gr_{\d} E \subseteq \c[x, \xi] \ .
\end{equation}
By \eqref{ass4}, all the codimensions here are finite; exactly 
as in the proof of Theorem~\ref{main}, it follows that $\, \D \,$ 
is Noetherian and that $\, E \,$ is finite over $\, \D \,$.  Thus 
all the assumptions of the Levasseur-Stafford lemma hold.

It remains to see that $\, V \in \Gr $.  We show first that $\, V \,$ 
is primary decomposable: according to Theorem~\ref{CH2},  
it is equivalent to show that 
$\, P.\c[x] = V \,$ (where $\, P \,$ is as in \eqref{P}).  
Now, if $\, D \in \D = \D(V) \,$, it follows from the definition 
\eqref{P} of $\, P \,$ that $\, DP \subseteq P \,$; hence 
$\, DP.\c[x] \subseteq P.\c[x] \,$, that is, 
$\, P.\c[x] \,$ is a left sub-$\D$-module of $\, V \,$.  
Next, let $\, p \,$ and $\, q \,$ be as in 
Proposition~\ref{grat}; then 
$\, qpV \subseteq q \c[x] \subseteq V \,$, so $\, qp \in \D \,$.  
Also, since $\, q \in P \,$, we have 
$$
qpV \subseteq Pp.V \subseteq P.\c[x] \ .
$$
Thus $\, qp \,$ is a nonzero element in the annihilator 
of the $\, \D$-module $\, V/{P.\c[x]} \,$.
Since $\, \D \,$ is simple, this annihilator must be all of $\, \D \,$, 
in particular it must contain $1 \,$.  This shows that 
$\, V/{P.\c[x]} \,$ is the zero module, that is, 
$\, P.\c[x] = V \,$, as claimed.

Finally, the fact that 
$\, V \in \Gr $ is a consequence 
of our assumption that the operators in $\, \Check{B} \,$ are 
normalized with first two 
coefficients constant. Indeed, since $\, V \,$ is primary decomposable, 
we have $\, V = fW \,$ for some $\, W \in \Gr, \, f \in \c(x) \,$. 
Clearly $\, \D(V) = f \D(W) f^{-1} \,$; thus we have 
$\, \D = f \D(W) f^{-1} \,$
Conjugating by $\, f \,$ 
does not change either the $\d$-  or the $x$-filtration on 
$\, \c(x)[\d_x] \,$, thus $\, \Check{B} =  f \mathcal{A}_{b(W)} f^{-1} \,$. 
So the algebras $\, \Check{B} \,$ and $\, f^{-1 }\Check{B} f \,$ both
consist of operators with first two coefficients constant.  But the 
second coefficients differ by nonzero multiples of $\, f' f^{-1} \,$, 
so that is possible only if $\, f \,$ is constant, and 
hence $\, V = W \,$.
\end{proof}

To end this section, we give the promised proof of 
the uniqueness of the space $\, V $ that we have constructed. 
The proof depends on the fact that  $\, \D \,$ is simple.
We note first
\begin{lemma}
\label{one}
There exists at most one (nonzero) simple sub-$\D$-module of $\, \c(x) \,$.
\end{lemma} 
\begin{proof}
Suppose that $\, V \,$ and $\, V' \,$ are two such modules. 
Since $\, \Frac B = \c(x) \,$, 
if we fix nonzero elements 
$\, v \in V, \, v' \in V' \,$, we can find nonzero $\, a, b \in B \,$ 
such that $\, av = bv' \,$.  Thus $\, V \cap V' \not= 0 \,$, 
hence (since $\, V \,$ and $\, V' \,$ are both simple) 
$\, V = V \cap V' = V' \,$. 
\end{proof}

The uniqueness of $\, V \,$ follows from Lemma~\ref{one} and 
the next Proposition.
\begin{proposition}
\label{last}
Let $\, V $ be any nonzero sub-$\D$-module of $\, \c(x) $ that is 
finite over $\, B \,$.  Then  $\, V \,$ is simple.
\end{proposition}
\begin{proof}
Let $\, U \subseteq V \,$ be a nonzero sub-$\D$-module: 
fix any nonzero element $\, u \in U \,$.  Let 
$\, \{ v_1, \ldots, v_k \} \,$ generate $\, V \,$ as a 
$B$-module.  Since $\, \Frac B = \c(x) \,$, we can find nonzero
$\, a_i, b_i \in B \,$ such that $\, a_i v_i = b_i u \,$ 
(for $\, 1 \leq i \leq k \,$).  It follows that if $\, a \,$ is the 
product of the $\, a_i \,$, then $\, a v_i \in U \text{ for all } i \,$, 
hence $\, aV \subseteq U \,$. Thus the annihilator (in $\, \D \,$) of 
$\, V/U \,$ is nonzero; because $\, \D \,$ is simple, it must be the 
whole of $\, \D \,$, so $\, U = V \,$.  Hence $\, V \,$ is simple.
\end{proof}

\section{Proof of Theorem~\ref{T4}}

For the rest of the paper $\, \D \,$ will denote an algebra which is 
isomorphic to  $\, \D(W) \,$ for some $\, W \in \Gr$. 
Let us define $\, \Grad \D \,$ to be the set of all isomorphisms 
$\, \sigma_W : \D(W) \to \D \,$ for various $\, W \in \Gr$ (more precisely, 
$\, \Grad \D \,$ is the set of all pairs $\, (W, \sigma_W) \,$ 
where $\, W \in \Gr$ and $\, \sigma_W \,$ is an isomorphism  
as above).  On the other 
hand, let $\, \Fad \D \,$ denote the set of all triples $\, (B,x,\d_x) \,$ 
where $\, B \,$ is a mad subalgebra of $\, \D \,$ and  $\, (x,\d_x) \,$
is a good fat framing of $\, B \,$.  
The set $\, \Fad \D(W) \,$ has the natural base-point 
$\, (A_W,z,\d_z) \,$: thus there is an obvious map  
$\, \alpha: \Grad \D \to \Fad \D \,$ 
which assigns to $\, \sigma_W \in \Grad \D \,$ the point 
$$
\alpha(\sigma_W) := (\sigma_W(A_W), \sigma_W(z), \sigma_W(\d_z)) \in \Fad \D 
$$
(the map $\, \sigma_W \,$ extends to an isomorphism of quotient fields 
$\, Q \to \Frac \D \,$, which we denote by the same symbol).
We can reformulate Theorem~\ref{main} as follows.
\begin{theorem}
\label{alpha}
The above map $\, \alpha : \Grad \D \to \Fad \D \,$ is bijective.
\end{theorem}
\begin{proof}
Let $\, (B,x,\d_x) \in \Fad \D \,$, and let 
$\, \theta : \c(x)[\d_x] \to \c(z)[\d_z] \,$ be the isomorphism 
which sends $\, x \,$ to $\, z \,$ and $\, \d_x \,$ to $\, \d_z \,$. 
Theorem~\ref{main} states that $\, \theta \,$ maps $\, \D \,$ 
isomorphically onto one of the algebras $\, \D(W) \,$, so the 
restriction of $\, \theta^{-1} \,$ to $\, \D(W) \,$ gives us 
a point of $\, \Grad \D \,$. It is 
clear that this construction defines the inverse map to $\, \alpha \,$.
\end{proof}

Observe now that the group $\, \Aut \D \times \Gamma \,$ acts 
naturally on each of spaces $\, \Grad \D \,$ and $\, \Fad \D \,$ 
(recall from the Introduction that $\, \Gamma \,$ is the group 
of maps $\, \gamma_p \,$ defined by \eqref{gamma}).  
Given $\, \sigma_W \in \Grad \D \,$, we can compose it with any 
$\, \sigma \in \Aut \D \,$ and $\, \gamma \in \Gamma \,$, as follows:
$$
\D(\gamma^{-1} W) \xrightarrow{\gamma} 
\D(W) \xrightarrow{\sigma_W} \D \xrightarrow{\sigma} \D \,,
$$
where the first map is explained in Subsection~\ref{gammact}.
This clearly defines an action of $\, \Aut \D \times \Gamma \,$ 
on $\, \Grad \D \,$.  The action of $\, \Aut \D \,$ on 
$\, \Fad \D \,$ is induced from its natural action on $\, \D \,$; 
we let $\, \gamma_p \in \Gamma \,$ act on $\, \Fad \D \,$ 
as formal conjugation by $\, e^{p(x)} \,$, that is, we set
\begin{equation}
\label{gamfad}
\gamma_p(B,x,\d_x) = (B,x,\d_x - p'(x)) \,.
\end{equation}
Directly from the definitions, we can check:
\begin{proposition}
\label{alpheq}
The bijection $\, \alpha \,$ in Theorem~\ref{alpha} is equivariant 
with respect to the above actions of $\, \Aut \D \times \Gamma \,$.
\end{proposition}

It follows that $\, \alpha \,$ induces bijections between the 
quotient spaces of $\, \Grad \D \,$ and $\, \Fad \D \,$ by any 
one of the groups $\, \Aut \D \,$, $\, \Gamma \,$ or 
$\, \Aut \D \times \Gamma \,$.  The latter two possibilities will 
yield Theorems~\ref{T3} and \ref{T4}, respectively, but we consider 
first the quotient by $\, \Aut \D \,$.  
The obvious map $\, \Grad \D \to \Gr$ 
(sending $\, \sigma_W \,$ to $\, W \,$) clearly induces an injection from 
$\, {\Grad \D}/{\Aut \D} \,$ into $\, \Gr $.  Its image consists of 
all $\, W \in \Gr $ such that $\, \D(W) \,$ is isomorphic to $\, \D \,$: 
as explained in \cite{BW1}, \cite{BW3}, this consists of one of the 
{\it Calogero-Moser strata} $\, \C_n \subset \Gr $.  We therefore 
obtain
\begin{corollary}
\label{cor1}
The bijection $\, \alpha \,$ of Theorem~\ref{alpha} induces a bijection 
$$
\C_n \to {\Fad \D}/{\Aut \D}  \ ,
$$
where $\, n \,$ is the integer determining the isomorphism class 
of $\, \D \,$.
\end{corollary}

We now divide out further by the action of $\, \Gamma \,$.  According to 
\cite{BW2}, the action of $\, \Gamma \,$ on the space $\, \C_n \,$ 
is as defined by \eqref{Caction}; while the formula \eqref{gamfad} shows 
that the quotient map $\, \Fad \D \to {\Fad \D}/{\Gamma} \,$ can be 
identified with the forgetful map from $\, \Fad \D \,$ to $\, \Mad \D \,$, 
sending $\, (B, x, \d_x) \,$ to $\, (B, x) \,$.  Hence 
Corollary~\ref{cor1} yields the following 
slightly sharpened version of Theorem~\ref{T4}.
\begin{corollary}
\label{gent4}
The bijection of Corollary~\ref{cor1} induces a bijection
$$
\C_n / \Gamma \to {\Mad \D} / {\Aut \D} \,.
$$
\end{corollary}

We can also divide out just by the action of $\, \Gamma \,$, giving 
a bijection from $\, {\Grad \D} / {\Gamma} \,$ to $\, \Mad \D \,$.  
As mentioned above, this leads to Theorem~\ref{T3}, but more 
work is needed to identify  $\, {\Grad \D} \,$ with the 
space $\, \Aut A \,$ in that theorem.

\section{Proof of Theorem~\ref{T3}}
\label{10}

To obtain Theorem~\ref{T3} from the considerations in the preceding section, 
we need one more ingredient; namely, we need to see that  
the obvious action of $\, \Aut \D \,$ on $\, \Grad \D \,$ extends to 
an action of the larger group  $\, \Pic \D \,$. We first review some 
general facts about $\, \Pic \D \,$ (which are valid for an 
arbitrary $\c$-algebra $\, \D \,$).  For more details, see 
(for example) \cite{B}, Chapter 2.

Recall that $\, \Pic \D \,$ is the group (under tensor 
product) of  isomorphism classes of invertible 
$\, \D$-$\D $-bimodules (over $\, \c \,$, that is, we consider only 
bimodules on which the left and right $\c$-vector space structures 
coincide). There is a natural homomorphism from $\, \Aut \D \,$ to 
$\, \Pic \D \,$ which assigns to $\, \sigma \in \Aut \D \,$ the bimodule
$\, {}_{\bar{\sigma}}\D_1 \,$ (that is, $\, \D \,$ itself, 
but with the left action twisted by the inverse $\, \bar{\sigma} \,$ 
of $\, \sigma \,$).  The kernel of 
this map is exactly the group of {\it inner} automorphisms of $\, \D \,$. 
At the cost of breaking the left/right symmetry, we can describe 
the elements of $\, \Pic \D \,$ in the following way.  Let $\, M \,$ 
be an invertible $\, \D$-$\D $-bimodule: if we momentarily forget 
the left action of $\, \D \,$, then $\, M \,$ becomes a (progenerative) 
right $\, \D$-module $\, M_{\D} \,$.  The forgotten left action of 
$\, \D \,$ is then defined by some isomorphism
$$
\sigma : \End_{\D}(M_{\D})  \to  \D \,,
$$
which is again unique up to composition with an inner 
automorphism of $\, \D \,$.  

Now we return to our case, where $\, \D \,$ is isomorphic to one of the 
algebras $\, \D(W) \,$.  In this case the remarks above can be simplified 
a little.  We note first
\begin{lemma}
\label{sp1}
The algebras $\, \D(W) \,$ (where $\, W \in \Gr $) have no nontrivial 
inner automorphisms.
\end{lemma}
\begin{proof}
A differential operator $\, D \in \D(W) \,$ can be invertible only 
if it has order zero, that is, if it is a function.  But by 
Proposition~\ref{grdw}, the only functions in $\, \D(W) \,$ are polynomials, 
hence the only invertible elements of $\, \D(W) \,$ are the scalars.
\end{proof}

It follows that we may regard $\, \Aut \D \,$ as a subgroup 
of $\, \Pic \D \,$ via the natural homomorphism described above.    
Next, we have the Cannings-Holland description of the ideal classes 
of $\, \D(W) \,$. 
\begin{lemma}
\label{grep}
For any $\, W \in \Gr $, 
each isomorphism class of right ideals of $\, \D(W) \,$ has a unique 
representative of the form
$$
\D(W,V) := \{ D \in \c(z)[\d_z] : D.W \subseteq V \} 
$$
with $\, V \in \Gr$.
\end{lemma}
\begin{proof}
In the case when $\, W = \c[z] \,$, so that 
$\, \D(W) \,$ is the Weyl algebra $\, A \,$, this is exactly 
Theorem~\ref{CH1}: each ideal class 
in $\, A \,$ has a unique representative of the form $ \, \D(\c[z],V) \,$.  
But $\, \D(W) \,$ is Morita equivalent to $\, A \,$ via the 
invertible bimodule $ \, \D(W, \c[z]) \,$; it follows that each ideal class 
in $\, \D(W) \,$ has a unique representative of the form
$$
\D(\c[z],V) \D(W, \c[z]) = \D(W,V) \, ,
$$
as claimed.
\end{proof}

With these preliminaries, we can define 
the action of $\, \Pic \D \,$ on $\, \Grad \D \,$. 
Let $\, [M] \in \Pic \D \,,\  \sigma_W \in \Grad \D \,$. It follows 
from Lemmas~\ref{sp1} and \ref{grep} that $\, [M] \,$ has a 
unique representative of the form $\, \D(W,V) \,$ with the structure 
of right $\, \D$-module determined via $\, \sigma_W \,$ and 
the structure of left $\, \D$-module determined via some isomorphism 
$\, \sigma_V : \D(V) \to \D \,$.  For short, in what follows we shall 
say that $\, [M] \,$ is {\it represented} by this triple 
$\, (\D(W,V), \sigma_W, \sigma_V) \,$. We define 
\begin{equation}
\label{action}
[M].\sigma_W = \sigma_V \,.
\end{equation}

\begin{theorem}
\label{picax}
The formula \eqref{action} defines a free transitive action of 
$\, \Pic \D \,$ on $\, \Grad \D \,$.
\end{theorem}
\begin{proof}
Straightforward.  The main point is to check that we do indeed 
have a group action, that is, if $\, [M] \, , [N] \in \Pic \D \,$ 
and $\, \sigma_W \in \Grad \D \, $, then 
$$
[N].([M].\sigma_W) = [N \otimes_{\D} M]. \sigma_W \,.
$$
That amounts to showing that if $\, [M] \,$ is represented by 
$\, (\D(W,V), \sigma_W, \sigma_V) \,$ and $\, [N] \,$ by 
$\, (\D(V,U), \sigma_V, \sigma_U) \,$ then $\, [N \otimes_{\D} M] \,$ 
is represented by $\, (\D(W,U), \sigma_W, \sigma_U) \,$. 
The map $\, D_1 \otimes D_2 \mapsto D_1 D_2 \,$ provides 
the required isomorphism of bimodules from $\, \D(V,U) \otimes \D(W,V) \,$ 
to $\, \D(W,U) \,$.  It is trivial to show that the action is free 
and transitive.
\end{proof}

Now recall that we have an action of the group $\, \Gamma \,$ 
on $\, \Grad \D \,$, commuting with the action of $\, \Aut \D \,$. 
A little more generally, we have

\begin{proposition}
\label{gamcom}
The above action of $\, \Pic \D \,$ on  $\, \Grad \D \,$ commutes 
with the action of $\, \Gamma \,$.
\end{proposition}
\begin{proof}
Let $\, [M] \in \Pic \D \,$ and $\, \sigma_W \in \Grad \D \,$. 
Let $\, [M].\sigma_W = \sigma_V \,$, so that
$\, [M] \,$ is represented by $\, (\D(W,V), \sigma_W, \sigma_V) \,$.  
If $\, \gamma \in \Gamma ,$ we have to 
show that $\, [M].(\sigma_W \gamma) = \sigma_V \gamma \,$; 
equivalently, that $\, [M] \,$ is also represented by 
$\, (\D(\gamma^{-1}W,\gamma^{-1}V), \sigma_W \gamma, \sigma_V \gamma) \,$. 
It is easy to check that the map
$$
\gamma : \D(\gamma^{-1}W,\gamma^{-1}V) = \gamma^{-1} \D(W,V) 
\to \D(W,V) 
$$
explained in Subsection~\ref{gammact} 
is a bimodule isomorphism; hence the result.
\end{proof}

Now let us fix a base-point 
$\, \sigma_W \in \Grad \D \,$; 
according to Theorem~\ref{picax}, the map
\begin{equation}
\label{picgr}
\Pic \D \to \Grad \D 
\end{equation}
which sends $\, [M] \,$ to $\, [M].{\sigma_W} \,$ is bijective.
Fixing a base-point gives us also 
a distinguished invertible $\, \D $-$ A $-bimodule 
$\, P := \D(\c[z],W) \,$, where it 
is understood that the structure of left $\, \D $-module on $\, P \,$ 
is defined via the isomorphism $\, \sigma_W \,$.  
By \eqref{uvw}, the inverse 
$\, A $-$ \D $-bimodule is $\, P^* := \D(W,\c[z]) \,$.  According to 
\cite{St}, the natural map $\, \Aut A \to \Pic A \,$ is 
an isomorphism; on the other hand, $\, P \,$ defines an isomorphism 
from $\, \Pic A \,$ to $\, \Pic \D \,$, sending (the class of) an 
$\, A $-$ A $-bimodule $\, M \,$ to   
$\, P \otimes_A M \otimes_A P^* \,$.  Combining the composite 
isomorphism $\, \Aut A \simeq \Pic \D \,$ with the bijection \eqref{picgr}, 
we obtain a bijective map
\begin{equation}
\label{autgr}
\beta : \Aut A \to \Grad \D \,.
\end{equation}
\begin{lemma}
\label{gcom}
Under the bijection $\, \beta \,$, the action of $\, \Gamma \,$ on 
$\, \Grad \D \,$ corresponds to its action by right multiplication on 
$\, \Aut A \,$.
\end{lemma}
\begin{proof}
Because of Proposition~\ref{gamcom}, it is enough to show that if 
$\, \gamma \in \Gamma \,$ then $\, \beta(\gamma) = \sigma_W \gamma \,$ 
(recall that $\, \sigma_W \,$ is our chosen base-point in 
$\, \Grad \D \,$).  Since $ \, \gamma \,$ corresponds to the bimodule 
$\, M := P \otimes_A {}_{\bar{\gamma}}A_1 \otimes_A P^* \,$ in 
$\, \Pic \D \,$, we have to see that this bimodule is represented by 
$\, (\D(W, \gamma^{-1}W), \sigma_W, \sigma_W \gamma) \,$.  It is easy to 
check that the map
$$
D_1 \otimes a \otimes D_2 \mapsto \gamma^{-1}(D_1) a D_2 
$$
defines the desired bimodule isomorphism from $\, M \,$ to 
$\, \D(W, \gamma^{-1}W) \,$.
\end{proof}

Finally, we can now consider the composite bijection
\begin{equation}
\label{ab}
\Aut A \xrightarrow{\beta} \Grad \D \xrightarrow{\alpha} \Fad \D \,.
\end{equation}
Using Lemma~\ref{gcom}, we can divide both sides 
of this bijection by $\, \Gamma \,$  
to obtain the following more precise version of 
Theorem~\ref{T3}.
\begin{theorem}
\label{final}
The bijection \eqref{ab} induces a bijection
$$
\Aut(A) / \Gamma \to \Mad \D \,.
$$
\end{theorem}
\begin{remark}
The bijection in Theorem~\ref{final} depends on the choice of base-point 
in $\, \Grad \D \,$; however, in practice there is usually a natural 
choice. For example, if $\, \D = \D(W) \,$ for some $\, W \in \Gr $, then 
it is natural to take the identity map in $\, \D(W) \,$ as the  base-point. 
Similarly, if $\, \D = \D(X) \,$, where 
$\, X \,$ is a framed curve with ring of functions 
$\, \O(X) \subseteq \c[z] \,$, 
then there is a unique monic polynomial $\, p(z) \,$ 
such that $\, W := p^{-1} \O(X) \,$ belongs to 
$\, \Gr $, and it is natural to take as 
base-point the isomorphism $\, \sigma_W : \D(W) \to \D(X) \,$ defined by 
$\, \sigma_W(D) := pDp^{-1} \,$.  This remark perhaps justifies our 
use of the word ``natural'' in Theorem~\ref{T3}.
\end{remark}

\section{Appendix: some properties of $\, \D(W) \,$}

Here we provide proofs that the algebras $\, \D(W) \,$ 
(for $\,W \in \Gr $) have the properties needed for us to apply the 
results of Sections~\ref{several} and \ref{dual}.  
The proofs of the first two propositions 
use only the fact that we have 
\begin{equation}
\label{pq}
p \c[z] \subseteq W \subseteq q^{-1} \c[z]
\end{equation}
for suitable polynomials $\, p,q \,$; thus these Propositions hold 
also for spaces $\, W \,$ in the larger Grassmannian $\, \Grat $. 
\begin{proposition}
\label{fracd}
The field of fractions of $\, \D(W) \,$ is $\, Q \,$.
\end{proposition}
\begin{proof}
It follows from \eqref{pq} that 
if $\, D \in \c[z, \d_z] \,$, then
$$
pDq.W \subseteq pD .\c[z] \subseteq p \c[z] \subseteq W \,,
$$
that is, $\, p \c[z, \d_z] q \subseteq \D(W) \,$.  It follows that 
the quotient field of $\, \D(W) \,$ contains the Weyl algebra 
$\, \c[z, \d_z] \,$; it is therefore the whole of $\, Q \,$, as claimed.
\end{proof}
\begin{proposition}
\label{grdw}
$ \gr_{\d} \D(W) \subseteq \c[z, \zeta] \,$.
\end{proposition}
\begin{proof}
Arguing as in the previous proof, it follows from \eqref{pq} that
$\, \D(W) \,$ is contained in $\, q^{-1} \c[z, \d_z] p^{-1} \,$; 
hence the leading coefficient of every element $\, L \in \D(W) \,$ 
has denominator at worst $\, pq \,$.  But this is true of 
$\, L^n \,$ for every $\, n \geq 1 \,$, hence that leading coefficient 
must be a polynomial, as claimed. 
\end{proof}

Our last Proposition (which is {\it not} valid for all $\, W \in \Grat $) 
is less easy to prove. The proof given in \cite{SS} for 
$\, \D(X) \,$ (where $\, X \,$ is a framed curve) generalizes easily 
to our case; here we give another proof, using the existence of the 
bispectral involution on $\, \Gr $.
\begin{proposition}
\label{finite}
The pair $\, (\D(W), A_W) \,$ satisfies the condition \eqref{ass4}.
\end{proposition}
\begin{proof}
We observed in the proof of Proposition~\ref{fracd} that
$\, p \c[z, \d_z] q \subseteq \D(W) \,$
for suitable polynomials $\, p \,$ and $\, q \,$; it follows that
\begin{equation}
\label{a}
p(z)q(z) \c[z, \zeta] \subseteq \gr_{\d} \D(W) \ .
\end{equation}
Similarly, there are polynomials $\, r \,$ and $\, s \,$ such that 
$\, r(z) \c[z, \d_z] s(z) \subseteq \D(b(W)) \,$, 
hence (applying the anti-automorphism $\, b \,$ of $\, Q $)
$$
s(\d_z) \c[z, \d_z] r(\d_z) \subseteq b \,\D(b(W)) = 
K_{b(W)}^{-1} \D(W) K_{b(W)} \ .
$$
Since $\, K-1 \,$ has negative $\d$-filtration, it follows that 
\begin{equation}
\label{c}
\zeta^N \c[z, \zeta] \subseteq \gr_{\d} \D(W) \ ,
\end{equation}
(where $\, N := \deg r + \deg s $). It follows at once from 
\eqref{a} and \eqref{c} that $\,   \gr_{\d} \D(W) \,$
has finite codimension in $\, \c[z, \zeta] \,$.
\end{proof}

\bibliographystyle{amsalpha}

\begin{thebibliography}{A}
%
\bibitem[AM]{AM}
M. F. Atiyah and I. G. Macdonald, 
\textit{Introduction to Commutative Algebra},
Addison-Wesley, 1969.
%
\bibitem[B]{B}
H. Bass, 
\textit{Algebraic K-Theory},
Benjamin, 1968.
%
\bibitem[BB]{BB}
A. Beilinson and J. Bernstein, \textit{A proof of Jantzen conjectures},
Advances in Soviet Mathematics \textbf{16} (1993), 1--50.
%
\bibitem[BC]{BC}
J. L. Burchnall and T. W. Chaundy, \textit{Commutative ordinary 
differential operators}, Proc. London Math. Soc. \textbf{21} 
(1923), 420--440.
%
%
\bibitem[BW1]{BW1}
Yu. Berest and G. Wilson, \textit{Classification of rings of 
differential operators on affine curves}, Internat. Math. Res. 
Notices 1999(2), 105--109.
%
\bibitem[BW2]{BW2}
Yu. Berest and G. Wilson, \textit{Automorphisms and ideals of the 
Weyl algebra}, Math. Ann. \textbf{318} (2000), 127--147.
%
\bibitem[BW3]{BW3}
Yu. Berest and G. Wilson, \textit{Differential isomorphism and 
equivalence of algebraic varieties}, in Proceedings of the 2002 
Oxford Symposium \textit{Topology, Geometry and Quantum Field Theory} 
(Ed.\ U.\ Tillmann), London Mathematical Society Lecture Notes 
Series \textbf{308}, Cambridge University Press, Cambridge,
2004, pp.\ 98--126.
%
\bibitem[CH]{CH}
R. C. Cannings and M. P. Holland, \textit{Right ideals of rings
of differential operators}, 
J. Algebra \textbf{167} (1994), 116--141.
%
\bibitem[D]{D}
J. Dixmier, \textit{Sur les alg\`ebres de Weyl}, Bull. Soc. Math. France 
\textbf{96} (1968), 209--242.
%
%
%
\bibitem[E]{E2}
P. Eakin, \textit{A note on finite-dimensional subrings of polynomial 
rings}, Proc. Amer. Math. Soc.
\textbf{31} (1972), 75--80.
%
%
%
\bibitem[J]{J}
N. Jacobson, \textit{Lectures in Abstract Algebra III}, 
Van Nostrand, 1964.
%
\bibitem[KM]{KM}
A. A. Klein and L. Makar-Limanov, \textit{Skew fields of differential 
operators}, Israel J. Math. \textbf{72} 1990, 281--287. 
%
\bibitem[K1]{K1}
K. M. Kouakou, \textit{Isomorphismes entre alg\`ebres d'op\'erateurs 
diff\'erentielles sur les courbes alg\'ebriques affines}, 
Th\`ese, Univ. Claude Bernard Lyon-I, 1994.
%
\bibitem[K2]{K2}
K. M. Kouakou, \textit{Classification des id\'eaux \`a droite de 
$A_1(k)$}, Bull. London Math. Soc. \textbf{35} (2003), 503--512.
%
%
\bibitem[LM]{LM}
G. Letzter and L. Makar-Limanov, 
\textit{Rings of differential operators over rational affine curves},
Bull. Soc. Math. France \textbf{118} (1990), 193--209.
%
\bibitem[LS]{LS}
T. Levasseur and J. T. Stafford, 
\textit{Rings of differential operators on classical rings of invariants}, 
Mem. Amer. Math. Soc. \textbf{412}, 1989.
%
%
\bibitem[M]{M}
L. Makar-Limanov, \textit{Rings of differential operators 
on algebraic curves}, Bull. London Math. Soc. \textbf{21} (1989),
538--540.
%
%
%
%
\bibitem[P]{P}
P. Perkins, \textit{Commutative subalgebras of the ring of 
differential operators on a curve}, Pacific J. Math.,
\textbf{139}(2) (1989), 279--302.
%
%
\bibitem[S]{S}
I. Schur, \textit{\"Uber vertauschbare lineare Differentialausdr\"ucke}, 
Sitzungsber. Berliner Math. Ges. \textbf{4} (1905), 2--8; 
Gesammelte Abhandlungen Band I, Berlin-Heidelberg-New York (1973), 
170--176.
%
\bibitem[SW]{SW}
G. Segal and G. Wilson,
\textit{Loop groups and equations of KdV type}, 
Publ. Math. IHES \textbf{61} (1985), 
5--65.
%
\bibitem[SS]{SS}
S. P. Smith and J. T. Stafford, \textit{Differential
operators on an affine curve}, Proc. London Math. Soc.
(3) \textbf{56} (1988), 229--259.
%
\bibitem[St]{St}
J. T. Stafford, \textit{Endomorphisms of right ideals of the Weyl algebra},
Trans. Amer. Math. Soc. \textbf{299} (1987), 623--639.
%
\bibitem[W1]{W1}
G. Wilson,
\textit{Bispectral commutative ordinary differential
operators}, J. reine angew. Math. \textbf{442} (1993), 
177--204.
%
\bibitem[W2]{W2}
G. Wilson,
\textit{Collisions of Calogero-Moser particles and
an adelic Grassmannian} (with an Appendix by I. G. Macdonald), 
Invent. Math. \textbf{133} 
(1998), 1--41.
%
\end{thebibliography}

\end{document}